\RequirePackage{fix-cm}
\documentclass[smallextended]{svjour3}       
\smartqed  
\usepackage{amsmath, amsfonts,amssymb, graphicx, cite}
\usepackage[titletoc,title]{appendix}

\newcommand{\e}{\varepsilon}

\newcommand{\x}{{\bf x}}
\newcommand{\z}{{\bf z}}
\newcommand{\n}{{\bf n}}
\newcommand{\R}{{\mathbb R}}

\newcommand{\bg}{{\bf g}}
\newcommand{\bh}{{\bf h}}
\newcommand{\ff}{{\bf f}}
\newcommand{\y}{{\bf y}}
\newcommand{\X}{{\bf X}}
\newcommand{\F}{{\bf F}}
\newcommand{\LL}{{\bf L}}
\newcommand{\T}{{\bf T}}
\newcommand{\B}{{\bf B}}

\newcommand{\yy}{\mathbf{y}}
\newcommand{\geucl}{{\bg_{\rm Eucl}}}
\newcommand{\dd}{{\rm d}}
\newcommand{\uu}{\mathbf{u}}
\newcommand{\vv}{\mathbf{v}}
\newcommand{\ww}{\mathbf{w}}
\newcommand{\sss}{\mathfrak{S}_{\rm CDS}}
\newcommand{\vvv}{\widetilde{\vv}}
\newcommand{\gggg}{\widetilde{\bg}}
\newcommand{\p}{\partial}
\newcommand{\na}{\nabla}

\def\be{\begin{equation}}
\def\ee{\end{equation}}
\def\bes{\begin{equation*}}
\def\ees{\end{equation*}}
\def\bc{\begin{cases}}
\def\ec{\end{cases}}

\numberwithin{equation}{section}
\numberwithin{theorem}{section}
\numberwithin{lemma}{section}
\numberwithin{proposition}{section}
\numberwithin{corollary}{section}
\numberwithin{definition}{section}
\numberwithin{remark}{section}
\numberwithin{example}{section}

\begin{document}

\title{Fluids, Elasticity, Geometry,\\ and the Existence of Wrinkled Solutions}
\author{Amit Acharya \and Gui-Qiang G. Chen \and Siran Li \and Marshall Slemrod \and \\ Dehua Wang}

\institute{A. Acharya \at
Civil \& Environmental Engineering, Carnegie Mellon University,
Pittsburgh, PA 15213, USA.\\
\email{acharyaamit@cmu.edu}
\and
G.-Q. Chen \at
Mathematical Institute, University of Oxford, Oxford, OX2 6GG, UK;\\
AMSS \& UCAS, Chinese Academy of Sciences, Beijing 100190, China.\\
\email{chengq@maths.ox.ac.uk}
\and
S. Li \at
Mathematical Institute, University of Oxford, Oxford, OX2 6GG, UK.\\
\email{siran.li@maths.ox.ac.uk}
\and
M. Slemrod \at
Department of Mathematics, University of Wisconsin, Madison, WI 53706, USA.\\
\email{slemrod@math.wisc.edu}
\and
D. Wang \at
Department of Mathematics, University of Pittsburgh, Pittsburgh, PA 15260, USA.\\
\email{dwang@math.pitt.edu}
}

\date{Received: May 3, 2016 / Accepted: May 15, 2017}
\maketitle

\begin{abstract} $\,\,$ We are concerned with$\,$ underlying connections$\,$ between fluids,
elasticity,
isometric embedding of Riemannian manifolds,
and the existence of wrinkled solutions of
the associated
nonlinear partial differential equations.
In this paper, we develop such connections for the case of two spatial dimensions,
and demonstrate that the continuum mechanical equations can be mapped into
a corresponding geometric framework and the
inherent direct application of the theory of isometric embeddings
and the Gauss-Codazzi equations
through examples for the Euler equations for fluids
and the Euler-Lagrange equations for elastic solids.
These results show that the geometric theory provides an avenue
for addressing the admissibility
criteria for nonlinear conservation laws in continuum mechanics.
\end{abstract}

\tableofcontents

\section{$\,$ Introduction}\label{S1}

We$\,$  are$\,$ concerned$\,$ with $\,$ underlying$\,$ connections$ \,$ between$\,$ fluids, $\,$ elasticity,
isometric embedding of Riemannian manifolds,
and the existence of wrinkled solutions of the associated
nonlinear partial differential equations.
One of the main purposes of this paper is to develop such connections
for the case of two spatial dimensions,
and examine whether the continuum mechanical equations can be mapped into
a corresponding geometric framework and the
inherent direct application of the theory of isometric embeddings and the Gauss-Codazzi equations.
Another motivation for such a study is to explore the possibility whether the geometric theory
can serve an avenue for addressing the admissibility criteria
for nonlinear conservation laws in continuum mechanics.

In recent years, a theory of {\it wild} solutions to
both the incompressible and compressible
Euler equations in two and higher spatial dimensions
has been developed
by De Lellis, Sz\'{e}kelyhidi Jr., and others
in
\cite{BDIS, BDS1, CDK, CDS2012, DS2009, DS2010, DS2012, DS2013, DS2014, DS2015,SzL2011,Wied}
and the references cited therein.
The approach is based on the
analogy with the highly irregular (static) solutions of the isometric
embedding of a two-dimensional Riemannian manifold into three-dimensional
Euclidean space given by the Nash-Kuiper theorem \cite{Kuiper,Nash1954}.
Specifically, the analogy arises because of the applicability of
Gromov's h-principle  and
convex integration to both problems ({\it cf}. \cite{Gromov1986}).
This suggests that perhaps the initial value problem in fluid dynamics
and the embedding problem in differential geometry would be more than just two
analogous issues, and could in fact be mapped one to the other.
A first thought
on this issue would suggest that the question is not even meaningful:
the fluid problem is dynamic and the embedding problem is static.
Thus, if the linkage is to make any sense at all,
we must think of the embedding problem as
dynamic and derive the equations of a time evolving two-dimensional
Riemannian manifold.
Of course, conceptually it is easy to visualize
this time evolving two-dimensional Riemannian manifold as a
two-dimensional surface moving in three-dimensional space
that can be seen in
such an everyday phenomenon as the vibration of the surface of a drum.

Thus, to make the link, we must derive the equations of a new type of geometric
flow and then interpret the consequences of this fluid-geometric duality. In
fact, once we have mapped a solution of the Euler equations onto the evolving
surface, it is rather easy to see that the dynamic metric $\bg$ is a short metric
in the sense of the Nash-Kuiper theorem \cite{Kuiper,Nash1954}
with respect to a metric associated with a developable surface.
Hence, an immediate consequence of our theory is that the evolving
manifold can be approximated by wrinkled manifolds
({\it i.e.},  manifolds with discontinuous second derivatives)
and thus gives some indication of the
existence of  {\it wild} solutions of the dual fluid problem.
A simple mental
picture of the geometric image of the fluid problem would be the motion of
$C^{1,\alpha}$ wrinkles on a piece of paper.
Unfortunately, while appealing, this mental picture is not correct, and
the correct visualization would be the fractal images given
by Borrelli et al. \cite{Borrelli2012,Borrelli2013}.
Moreover, the {\it wild} solutions of the geometric problem are
completely time-reversible and are in
reality a sequence of the Nash-Kuiper solutions of the embedding problem
with the metric given by $\bg=\bg^*$,
where $\bg^*$ corresponds to a developable surface
that is time-independent and has non-vanishing mean curvature.

The implication of this fact is immediate:
The geometric Nash-Kuiper {\it wild} solutions are time-reversible.
This suggests a plausible answer to the question raised
in \cite{BDIS, BDS1, CDK, CDS2012, DS2009, DS2010, DS2012, DS2013, DS2014, DS2015,SzL2011,Wied}
as to which is the correct admissibility criterion to choose the relevant solution
from the infinite number of non-unique solutions to the Euler equations.
Namely, no
dynamic admissibility criterion such as the energy inequality, entropy
inequality, entropy rate can serve the purpose,
since all the inequalities
become the identities for such solutions.
The only possible useful criteria
must be meaningful for time
reversible
fluid flow, such as  energy
minimization, artificial viscosity \cite{Dafermos-book},
or viscosity-capillarity \cite{Slemrod2013}.

The equations for geometric flow are abundant:
The Einstein equations of general
relativity and the Ricci flow equations are two of the better known examples.
In both of
these problems,
the metric for a Riemannian manifold becomes the dynamic
unknown.
In the theory developed in this paper, the same situation
arises: A dynamic metric $\bg$ is our unknown along with the second
fundamental form for the evolving manifold.
Furthermore, just as in the case
of the Einstein equations ({\it cf.} \cite{BIsen,KlaiR}),
the initial data must be
consistent with the problem, that is,
initially
an embedded
manifold does exist and the map is also consistent with the divergence free
condition on the velocity for the incompressible fluid case.
For the Einstein equations, this consistency of the initial data
yields the Einstein
constraint equations, while in our case  a system of constraint
equations is also required.
In this paper we show that our constraint equations have
local analytic solutions; moreover, there is a velocity $(u,v)$
that defines both a solution of the incompressible Euler equations and an
evolving two-dimensional Riemannian manifold isometrically immersed
in $\mathbb{R}^{3}$.
A similar result is given for the compressible case,
as well as for neo-Hookean elasticity.

The most important issue is the physical meaning of the
results mentioned above.
In short, we emphasize the
comments which have appeared in Sz\'{e}kelyhidi Jr. \cite{SzL2011},
Bardos-Titi-Wiedemann \cite{BTitiW}, Bardos-Titi \cite{BTiti},
and Bardos-Lopes Filho-Niu-Nussenzveig Lopes-Titi \cite{BLFN}.
The appearance of the {\it wild} solutions is due essentially
to the application of the Euler equations with vortex sheet initial data.
Hence, the Euler equations, which have no small
scales built into their theory as opposed to the compressible or
incompressible Navier-Stokes equations, have been used in a case
where they
should not be expected to apply.
To this, we must add the
proviso that was again mentioned in
De Lellis-Sz\'{e}kelyhidi Jr. \cite{DS2009, DS2010, DS2012, DS2013, DS2014, DS2015,SzL2011};
also see Elling \cite{Elling2006,Elling2013}:
these {\it wild} solutions could be a demonstration of fluid turbulence.
The results given in this paper show that, if this is the case, then this
{\it wild} fluid turbulent behavior is mirrored by the solutions generated by the
non-smooth moving wrinkled surfaces produced by the Nash-Kuiper theorem \cite{Kuiper,Nash1954}.

This paper consists of nine sections after this brief introduction.
In \S 2, we recall the Euler equations for an inviscid incompressible fluid and
the Gauss-Codazzi equations for the isometric embedding of a two-dimensional
Riemannian manifold into three-dimensional Euclidean space.
In \S 3, as discussed above, we derive the equations of the geometric flow,
as well as the constraint conditions on the initial data.
In \S 4, we prove the solvability of the constraint equations for the
initial data.  We emphasize that the constraint conditions require the
initial fluid velocity \textit{not} to be a shear flow. This assumption pairs
nicely with a theorem in De Lellis-Sz\'{e}kelyhidi Jr. \cite{DS2012} that non-smooth
shear flow initial data for both the compressible and incompressible
Euler equations are {\it wild} data and yield non-unique solutions to the Cauchy
problem for both the compressible and incompressible Euler equations.
In fact, our computation suggests that this is the only {\it wild} initial
data. On the other hand, we note that, for the degenerate case, when the fluid
motion is a shear flow, we still have a metric that provides the desired map,
namely, metric $\bg^*$.
In \S 5, we state and prove our main result:
Evolving from the initial data, there exists  a metric $\bg$
of the geometric flow equations which yields a solution of both the Euler equations
and the equations describing the evolving isometrically immersed Riemannian
manifold $({\mathcal M},\bg)$.
In \S 6, we  continue the discussion of the initial
data issue and show that, for the shear flow initial data in the hypotheses of
Lemma \ref{L41}, the symbol of the underlying system of second order partial
differential equations vanishes so that the equations for $\bg$
are degenerate;
however, as we just commented, metric $\bg^*$ suffices
in this singular case.
In \S 7, we give our principal result:
The evolving manifold arising from the Euler equations can be approximated by
wrinkled $C^{1,\alpha }$ manifolds $({\mathcal M}, \bg^*)$ for some $\alpha\in (0,1)$
which are continuous in time.
Since the time-continuity follows from a rather lengthy argument,
its proof is presented in a separate appendix of this paper for completeness.
These wrinkled solutions can be arranged as a time-sequence of solutions
which render the initial value problem for the evolving manifold to have an
infinite number of constant energy solutions.
Furthermore, we provide a formal map from the geometric wrinkled solutions to
weak solutions of the incompressible Euler equations.
In \S 8, a short discussion is provided to show that many of
our earlier results for the incompressible
Euler equations carry over to the compressible case.
Based on the knowledge that has been obtained from the fluid equations,
we show in \S 9 how the case of general continuum mechanics can be placed in our mechanics-geometry
framework. As an illustrative example, we demonstrate results for elastodynamic motion
of a neo-Hookean solid. Furthermore, these results
suggest that more refined continuum mechanical theories relying on micro-structure
could play a key role in choosing admissible solutions.
Our last section, \S 10, provides a discussion of
admissibility criteria for the (incompressible and compressible) Euler
equations which asserts that, for the multidimensional Euler equations, no dynamic admissibility condition would
eliminate {\it wild} solutions, and hence the only meaningful one must be the one which
eliminates or at least reduces the number of wrinkles of our dual geometric
problem. We also show how our work suggests a minimal dynamical model for internally stressed elastic materials.
The paper concludes with an appendix in which the time-continuity of the wrinkled
solutions is proved.

\section{$\,$ Basic Equations}\label{S2}

\subsection{$\,$ Geometric equations and notations}

We start with some basic geometric equations and notations for subsequent developments.
For more details, see Han-Hong \cite{HanHong} and the references cited therein.

Let $({\mathcal M},\bg)$ be a two-dimensional Riemannian manifold with
$\y(x_{1},x_{2})\in\R^{3}$ denoting a point on the manifold,
$\partial_{i}\y\cdot \partial_{j}\y=g_{ij}$.
The unit normal vector $\n$ to the manifold is given by
$$
\n=\frac{\partial_{i}\y\times \partial_{j}\y}{|\partial_{i}\y\times\partial_{j}\y|},
$$
and the second fundamental form is
$$
II=L(dx_{1})^{2}+2Mdx_{1}dx_{2}+N(dx_{2})^{2},
$$
where $L=\n\cdot \partial_{11}\y, \ M=\n\cdot \partial_{12}\y$,  $N=\n\cdot \partial_{22}\y$,
and  $\partial_{ij}:=\partial_i\partial_j$ with $\partial_i=\partial_{x_i}$
for $i,j=1,2$.
We will use the alternative version of the second fundamental form
$$
(L^{\prime}, M^{\prime}, N^{\prime})=\frac{1}{\sqrt{\det \bg}}(L,M,N),
$$
and recall that
\begin{equation}\label{2.1aa}
R_{1212}=\frac{\kappa}{\det \bg},
\end{equation}
where $\kappa$ is the Gauss curvature, and
$R_{ijkl}$ is the Riemann curvature tensor.
Furthermore, for notational simplicity, we henceforth drop the ``\'{}" superscript
in the alternative version of the second fundamental form.

Then the Codazzi equations are
\begin{equation}\label{e21}
\begin{cases}
\partial _{1}N-\partial_{2}M
=-\Gamma _{22}^{1}L+2\Gamma_{12}^{1}M-\Gamma_{11}^{1}N,\\[2mm]
\partial_{1}M-\partial_{2}L
 =\Gamma_{22}^{2}L-2\Gamma_{12}^{2}M+\Gamma_{11}^{2}N,
\end{cases}
\ee
and the Gauss equation is
\be\label{e22}
LN-M^{2}=\kappa.
\ee
The Christoffel symbols are given by the following formulas:
\begin{equation*}
\Gamma_{ij}^{k}={\frac12}g^{kl}(\partial_{j}g_{il}+\partial_{i}g_{jl}
-\partial_{l}g_{ij}),
\end{equation*}
so that
\be\label{e23}
\begin{cases}
\Gamma_{11}^{1}
=\frac1{2\det \bg}\left[g_{22}(\partial_{1}g_{11})-g_{12}(2\partial_{1}g_{12}-\partial_{2}g_{11})\right],\\[1mm]
\Gamma_{12}^{1}
=\frac1{2\det \bg}\left[g_{22}(\partial_{2}g_{11})-g_{12}(\partial_{1}g_{22})\right],\\[1mm]
\Gamma_{22}^{1}
=\frac1{2\det \bg}\left[g_{22}(2\partial_{2}g_{21}-\partial_{1}g_{22})-g_{12}(\partial_{2}g_{22})\right],\\[1mm]
\Gamma_{11}^{2}
=\frac1{2\det \bg}\left[-g_{12}(\partial_{1}g_{11})+g_{11}(2\partial_{1}g_{12}-\partial_{2}g_{11})\right],\\[1mm]
\Gamma_{12}^{2}
=\frac1{2\det \bg}\left[-g_{12}(\partial_{2}g_{11})+g_{11}(\partial_{1}g_{22})\right],\\[1mm]
\Gamma_{22}^{2}
=\frac1{2\det \bg}\left[-g_{12}(2\partial_{2}g_{21}-\partial_{1}g_{22})+g_{11}(\partial_{2}g_{22})\right],
\end{cases}
\ee
and the Riemann curvature tensor is given by
\begin{equation}\label{2.4-a}
R_{iljk}=g_{lm}\big(\partial_{k}\Gamma_{ij}^{m}
-\partial_{j}\Gamma_{ik}^{m}+\Gamma_{ij}^{n}\Gamma_{nk}^{m}
-\Gamma_{ik}^{n}\Gamma_{nj}^{m}\big),
\end{equation}
where we have used the Einstein summation convention that the repeated indices
are implicitly summed over in the terms, which will also be used from now on.

In particular, we have
\begin{equation}\label{2.4-b}
\begin{split}
R_{1212}
=&g_{21}\big(\partial_{2}\Gamma_{11}^{1}-\partial_{1}\Gamma_{12}^{1}
+\Gamma_{11}^{1}\Gamma_{12}^{1}+\Gamma_{11}^{2}\Gamma_{22}^{1}
-\Gamma_{12}^{1}\Gamma_{11}^{1}-\Gamma_{12}^{2}\Gamma_{21}^{1}\big)\\
&+g_{22}\big(\partial_{2}\Gamma_{11}^{2}-\partial_{1}\Gamma_{12}^{2}
+\Gamma_{11}^{1}\Gamma_{12}^{2}+\Gamma_{11}^{2}\Gamma_{22}^{2}
-\Gamma_{12}^{1}\Gamma_{11}^{2}-\Gamma_{12}^{2}\Gamma_{21}^{2}\big).
\end{split}
\end{equation}

With \eqref{2.1aa},
we have Gauss's Theorema
Egregium for the Gauss curvature.
A convenient form is given by Brioschi's
formula:
\be\label{e24}
\begin{split}
\kappa =&\frac1{2(\det \bg)^{2}}\det
\begin{bmatrix}
-
\partial_{22}g_{11}+2\partial_{12}g_{12}-
\partial_{11}g_{22} &
\partial_{1}g_{11}
& 2\partial_{1}g_{12}-
\partial_{2}g_{11} \\[1mm]
2\partial_{2}g_{12}-%
\partial_{1}g_{22}
& 2g_{11} & 2g_{12} \\[1mm]
\partial_{2}g_{22} & 2g_{12} & 2g_{22}%
\end{bmatrix} \\[2mm]
&-\frac1{2(\det \bg)^{2}}\det
\begin{bmatrix}
0 &
\partial_{2}g_{11} &
\partial_{1}g_{22} \\[1mm]
\partial_{2}g_{11} & 2g_{11} & 2g_{12} \\[1mm]
\partial_{1}g_{22} & 2g_{12} & 2g_{22}%
\end{bmatrix}.
\end{split}
\ee
We recall the fundamental theorem of surface theory states
that the solvability of the Gauss-Codazzi equations is a necessary and
sufficient condition for the existence of an isometric embedding,
{\it i.e.}, a
simply connected surface $\y\in\mathbb{R}^{3}$
which satisfies $\partial_{i}\y\cdot \partial_{j}\y=g_{ij}$.
A convenient reference for the smooth version of the fundamental theorem
is do Carmo \cite{CarmoM}, while a non-smooth version can be found
in
Mardare \cite{Mardare1, Mardare2}.

\subsection{$\,$ Incompressible Euler equations for an inviscid fluid}

The equations for the balance of linear momentum are
\be\label{e25}
\begin{cases}
\partial_{1}(u^{2}+p)+\partial_{2}(uv)=-\partial_{t}u,\\
\partial_{1}(uv)+\partial_{2}(v^{2}+p)=-\partial_{t}v,
\end{cases}
\ee
where the constant density $\rho =1$ is taken.
The condition of incompressibility for the constant density is then given by
the equation:
\be\label{e26}
\partial_{1}u+\partial_{2}v=0.
\ee
In addition, taking the divergence of equations \eqref{e25}
and using the incompressibility condition \eqref{e26},
we have
\be\label{e27}
\partial_{11}(u^{2})
+2\partial_{12}(uv)+\partial_{22}(v^{2})
=-\triangle p.
\ee

\subsection{$\,$ The geometric equations in fluid variables}

Just as in Chen-Slemrod-Wang \cite{CSW2010},
it is convenient to write the geometric
equations in fluid variables. Set
\be\label{e28}
L=v^{2}+p,\quad M=-uv, \quad N=u^{2}+p.
\ee
Then the Gauss equation becomes
$$
(v^{2}+p)(u^{2}+p)-(uv)^{2}=\kappa,
$$
that is,
$$
p^{2}+pq^{2}=\kappa,
$$
where $q^{2}=u^{2}+v^{2}$.
This quadratic equation then gives
\be\label{e29}
p=-\frac12q^{2}\pm \frac12\sqrt{q^{4}+4\kappa}.
\ee
This means that
$$
q^{4}+4\kappa \geq 0
$$
must be required.
As it will be seen in the analysis below,
this condition is always satisfied.

We have just shown that $(L,M,N)$ can be written in fluid variables.
Now we can write\ the fluid variables in terms of the geometric variables
$(L,M,N)$.
To do this, simply substitute formula \eqref{e29}
into \eqref{e28}
to find
$$
L-N=v^{2}-u^{2}.
$$
Write $v=-\frac{M}{u}$ to see that
$
L-N=\big(\frac{M}{u}\big)^{2}-u^{2},
$
which yields that
$$
u^{4}+(L-N)u^{2}-M^{2}=0.
$$
Then we see
$$
u^{2}=\frac12\left(-(L-N)\pm \sqrt{(L-N)^{2}+4M^{2}}\right),
$$
and, using $L-N=v^{2}-\big(\frac{M}{v}\big)^{2}$,
$$
v^{2}=\frac12\left((L-N)\pm \sqrt{(L-N)^{2}+4M^{2}}\right).
$$
This shows that the ``+'' sign must be chosen in the above formulas so that
\be\label{e210}
\begin{cases}
u^{2}=\frac12\left(-(L-N)+\sqrt{(L-N)^{2}+4M^{2}}\right),\\[1.5mm]
v^{2}=\frac12\left((L-N)+\sqrt{(L-N)^{2}+4M^{2}}\right).
\end{cases}
\ee
Note that the condition: $q^{4}+4\kappa \geq 0$, with $q^{2}=\sqrt{(L-N)^{2}+4M^{2}}$, is equivalent to
$$
\big((L-N)^{2}+4M^{2}\big)+4\kappa \geq 0,
$$
that is,
$$
\big((L-N)^{2}+4(LN-\kappa)\big)+4\kappa =(L+N)^{2}\geq 0,
$$
which is always satisfied.
Thus, we have shown that $(u,v,p)$ are determined by $(L,M,N)$,
since
$$
p=\frac{1}{2}\big(-q^{2}\pm \sqrt{q^{4}+4\kappa}\big),
\quad
q^{2}=u^{2}+v^{2}=\sqrt{(L-N)^{2}+4M^{2}}.
$$

\section{$\,$ The Equations for Geometric Flow}\label{S3}

We now construct a {\it dual} solution, which simultaneously satisfies
the incompressible Euler equations and the Gauss-Codazzi equations of
isometric embeddings.

As before, the Gauss-Codazzi equations are
\be\label{e31}
\begin{cases}
\partial_{1}N-\partial_{2}M
 =-\Gamma_{22}^{1}L+2\Gamma_{12}^{1}M-\Gamma_{11}^{1}N,\\[1mm]
\partial_{1}M-\partial_{2}L=\Gamma_{22}^{2}L-2\Gamma_{12}^{2}M+\Gamma_{11}^{2}N,
\end{cases}
\ee
and
\be\label{e32}
LN-M^2=\kappa.
\ee
Hence,  for any $(u,v)$
to be a solution of the Euler
equations \eqref{e25}--\eqref{e26}, we must have
\be\label{e33}
\begin{cases}
\partial_{t}u=\Gamma_{22}^{1}L-2\Gamma_{12}^{1}M+\Gamma_{11}^{1}N,\\[1mm]
\partial_{t}v=\Gamma_{22}^{2}L-2\Gamma_{12}^{2}M+\Gamma_{11}^{2}N,
\end{cases}
\ee
and
\be\label{e34}
\partial_{1}u+\partial_{2}v=0.
\ee
Taking the divergence of \eqref{e33} and using \eqref{e34}, we have
\be\label{e35}
\partial_{1}\big(\Gamma_{22}^{1}L-2\Gamma_{12}^{1}M+\Gamma_{11}^{1}N\big)
+\partial_{2}\big(\Gamma_{22}^{2}L-2\Gamma_{12}^{2}M+\Gamma_{11}^{2}N\big)=0.
\ee
For convenience, define
\be\label{e36}
\begin{split}
&U(L,M,N):=\left(\frac12\left[-(L-N)+\sqrt{(L-N)^{2}+4M^{2}}\right]\right)^{{\frac12}},\\
&V(L,M,N):=\left(\frac12\left[(L-N)+\sqrt{(L-N)^{2}+4M^{2}}\right]\right)^{{\frac12}},
\end{split}
\ee
so that
\be\label{e37}
(u, v)=(U(L,M,N),V(L,M,N)).
\ee
The other choice of
$$
(u, v)=-(U(L,M,N),V(L,M,N))
$$
can be handled similarly.

\smallskip
In summary, we have the two evolution equations \eqref{e33} and four closure
relations \eqref{e31}--\eqref{e32} and \eqref{e35}.
Since $\Gamma _{ij}^{k}$ and $\kappa$ are the functions of $(L, M, N, \bg)$ through \eqref{e23}--\eqref{e24},
we obtain the six equations \eqref{e31}--\eqref{e33} and \eqref{e35} for the six unknowns $(L, M, N, \bg)$.

\section{$\,$ The Constraint Equations and Their Consequences}\label{S4}

In this section, we exposit the constraint equations on the initial data and
the consequences of their solvability.
Our first result is

\begin{theorem}\label{Theorem41}
Assume that the initial data $(L, M, N, \bg)$ at $t=0$ satisfy the five constraint equations:
\begin{align}
&\partial_{1}U+\partial_{2}V=0, \label{e41}\\
&\partial_{1}\big(\Gamma_{22}^{1}L-2\Gamma_{12}^{1}M+\Gamma_{11}^{1}N\big)
  +\partial_{2}\big(\Gamma_{22}^{2}L-2\Gamma_{12}^{2}M+\Gamma_{11}^{2}N\big)=0,\label{e42}\\
&\partial_{1}N-\partial_{2}M=-\Gamma_{22}^{1}L+2\Gamma_{12}^{1}M-\Gamma_{11}^{1}N, \label{e43}\\
&\partial_{1}M-\partial_{2}L=\Gamma_{22}^{2}L-2\Gamma_{12}^{2}M+\Gamma_{11}^{2}N,\label{e44}\\
&LN-M^{2}=\kappa, \label{e45}
\end{align}
which mean that the initial data are consistent with the
incompressibility of the fluid and that the evolving manifold is initially
indeed a Riemannian manifold.
Then, if the system of six equations \eqref{e31}--\eqref{e33} and \eqref{e35}
in the six unknowns \\
$(L, M, N, \bg)$ is
satisfied, it produces simultaneously a solution of both the Gauss-Codazzi
equations and the incompressible Euler equations.
\end{theorem}

\proof
$\,$ Since $(u,v)=(U(L,M,N),V(L,M,N))$,  we see from the first two
evolution equations \eqref{e31} and \eqref{e35} that
$$
\partial_{t}(\partial_{1}U+\partial_{2}V)=0.
$$
Since the initial data satisfy the constraint equations,
we conclude that, for $t>0$,
$$
\partial_{1}U+\partial_{2}V=0.
$$

Next, since the Gauss-Codazzi equations are satisfied for all $t>0$,
the fluid representation \eqref{e28} for $(L, M, N)$
allows
to write the Codazzi equations as
\begin{eqnarray*}
&\partial_{1}(u^{2}+p)+\partial_{2}(uv)
  =-\Gamma_{22}^{1}L+2\Gamma_{12}^{1}M-\Gamma_{11}^{1}N,\\
&\partial_{1}(uv)+\partial_{2}(v^{2}+p)
  =-\Gamma_{22}^{2}L+2\Gamma_{12}^{2}M-\Gamma_{11}^{2}N,
\end{eqnarray*}
and the Gauss equation as
$$
p=\frac{1}{2}\big(-q^{2}+\sqrt{q^{4}+4\kappa}\big).
$$
By \eqref{e33},
we see that the balance of linear momentum equations is also satisfied.
The proof is complete.
\endproof

We next examine the solvability of the initial data system.
Note that \eqref{e42}--\eqref{e44} imply
\be\label{e46}
\partial_{1}(\partial_{1}N-\partial_{2}M)
+\partial_{2}(-\partial_{1}M+\partial_{2}L)=0,
\ee
and the initial data system can be written as
\begin{align}
&\partial_{1}U+\partial_{2}V=0, \label{e41b}\\
&\partial_{11}N-2\partial_{12}M+\partial_{22}L=0,\label{e46b}\\
&\partial_{1}N-\partial_{2}M=-\Gamma_{22}^{1}L+2\Gamma_{12}^{1}M
        -\Gamma_{11}^{1}N, \label{e43b}\\
&\partial_{1}M-\partial_{2}L=\Gamma_{22}^{2}L-2\Gamma_{12}^{2}M
   +\Gamma_{11}^{2}N,\label{e44b}\\
&LN-M^{2}=\kappa. \label{e45b}
\end{align}
We can reverse the above computation. If \eqref{e46b}--\eqref{e44b} are satisfied,
we
take the divergence of the left-hand sides of \eqref{e43b}--\eqref{e44b}
and employ \eqref{e46b}
to yield \eqref{e42}.

From \eqref{e23},
we find that
\eqref{e41} and \eqref{e43}--\eqref{e46} becomes an
undetermined
system of five equations in the six unknowns $(L, M, N, \bg)$.

Our existence result for the initial data satisfying the constraint equations
reads as follows:

\begin{lemma}\label{L41}
Let an analytic divergence free velocity $(u, v)$ be
prescribed in a neighborhood of a point $(x_{1,}x_{2})=(0,0)$
such that $uv\neq 0$ at this point.
Set $(x_{1}^{\prime}, x_{2}^{\prime})=(x_{1}+x_{2}, x_{1}-x_{2})$.
On $x_{1}^{\prime}=0$ {\rm (}respectively $x_{2}^{\prime}=0${\rm )},
prescribe the analytic initial data{\rm :}
$g_{11}=g_{22}=1,
\partial_{x_{1}^{\prime}}g_{11}=\partial_{x_{1}^{\prime}}g_{22}=0$,
and $g_{12}$ satisfying the ordinary
differential equation:
\begin{equation}\label{ode-1}
\begin{split}
\partial_{x_{2}^{\prime}}g_{12}=
&\frac1{2N}\big[g_{12} (\partial_{1}N-\partial_{2}M)+ (\partial_{1}M+\partial_{2}L)\big]\\
&-\frac1{2L}\big[(\partial_{1}N-\partial_{2}M)+g_{12} (\partial_{1}M+\partial_{2}L)\big],
\end{split}
\end{equation}
with initial condition $g_{12}=0$ at $x_{2}^{\prime }=0$
so that $g_{ij}=\delta _{ij}$ at $(x_{1}^{\prime}, x_{2}^{\prime})=(0,0)$
{\rm (}respectively, $\partial_{x_{2}^{\prime}}g_{11}=\partial_{x_{2}^{\prime}}g_{22}=0$
and a similar ordinary differential equation
and initial data{\rm )}.
Then the initial data system \eqref{e41}--\eqref{e45}
has a local analytic solution $(g_{11}, g_{12}, g_{22})$.
\end{lemma}

\proof $\,$  We divide the proof into eight steps.

\smallskip
{\bf 1.} In fluid variables, three of our equations are
\begin{align}
&\partial_{1}u+\partial _{2}v=0, \label{e41f} \\
&\partial_{11}(u^{2})+2\partial_{12}(uv)
+\partial_{22}(v^{2})=-\triangle p, \label{e46f}\\
&p=\frac12\big(-q^{2}+\sqrt{q^{4}+4\kappa}\big). \label{e45f}
\end{align}
Next, prescribe the velocity to make this a determined system. If a
divergence free velocity $(u, v)$ is prescribed,
then \eqref{e46f}  is immediately
solvable for $p$ under the standard regularity assumptions
on $(u, v)$, and hence \eqref{e45f}
defines the Gauss curvature $\kappa$.
Thus, $(L, M, N)$ are known to be independent of metric $\bg$.

\smallskip
{\bf 2.} Now, determine $\bg$ by solving the three equations:
\begin{align}
&\kappa \det \bg=g_{21}\big(\partial_{2}\Gamma_{11}^{1}
    -\partial_{1}\Gamma_{12}^{1}+\Gamma _{11}^{1}\Gamma_{12}^{1}
    +\Gamma_{11}^{2}\Gamma_{22}^{1}-\Gamma_{12}^{1}\Gamma_{11}^{1}
     -\Gamma _{12}^{2}\Gamma_{21}^{1}\big)\notag \\
&\qquad \qquad + g_{22}\big(\partial_{2}\Gamma_{11}^{2}
  -\partial_{1}\Gamma_{12}^{2}+\Gamma_{11}^{1}\Gamma_{12}^{2}
  +\Gamma_{11}^{2}\Gamma_{22}^{2}-\Gamma_{12}^{1}\Gamma_{11}^{2}
   -\Gamma_{12}^{2}\Gamma_{21}^{2}\big), \label{e48}\\
&\partial_{1}N-\partial_{2}M=-\Gamma_{22}^{1}L+2\Gamma_{12}^{1}M
   -\Gamma_{11}^{1}N, \label{e43c} \\
&-\partial_{1}M+\partial_{2}L=-\Gamma_{22}^{2}L+2\Gamma_{12}^{2}M
  -\Gamma_{11}^{2}N. \label{e44c}
\end{align}

\smallskip
{\bf 3.} We use \eqref{e43c}--\eqref{e44c} to
solve for $\partial_{1}g_{12}$ and $\partial_{2}g_{21}$.
Simply use  \eqref{e23}
to write these two equations \eqref{e43c}--\eqref{e44c} as
\begin{align}
&\frac1{\det \bg}
\big[
-g_{22}L (\partial_{2}g_{21})
+
g_{12}N(\partial_{1}g_{12})\big]\notag\\
&=\partial_{1}N-\partial_{2}M
+\frac1{2\det \bg}\left[
g_{22}(-\partial _{1}g_{22})-%
g_{12}(\partial _{2}g_{22})\right]L\notag\\
&\quad -\frac{1}{\det \bg}\left[%
g_{22}(\partial_{2}g_{11})-%
g_{12}(\partial _{1}g_{22})\right]M\notag\\
&\quad +\frac1{2\det \bg} \left[
g_{22}(\partial _{1}g_{11})-%
g_{12}(-\partial _{2}g_{11})\right] N, \label{e421}
\end{align}
\begin{align}
&\frac1{\det \bg}
\big[
g_{12}L(\partial_{2}g_{21})
-
g_{11}N(\partial_{1}g_{12})\big]\notag\\
&=-\partial_{1}M+\partial_{2}L
+\frac1{2\det \bg}\left[
g_{12}(-\partial_{1}g_{22})+%
g_{11}(\partial _{2}g_{22})\right]L\notag\\
&\quad-\frac{1}{\det \bg}
\left[ -
g_{12}(\partial_{2}g_{11})+%
g_{11}(\partial _{1}g_{22})\right] M\notag\\
&\quad+\frac1{2\det \bg}\left[ -%
g_{12}(\partial _{1}g_{11})+%
g_{11}(-\partial _{2}g_{11})\right]N.\label{e422}
\end{align}
In the matrix form, we have
\bes
\begin{split}
\frac1{\det \bg}
\begin{bmatrix}
g_{12}N  & -g_{22}L  \\
-g_{11}N  & g_{12}L
\end{bmatrix}%
\begin{bmatrix}
\partial_{1}g_{12} \\
\partial_{2}g_{21}%
\end{bmatrix}
= \begin{bmatrix}
\partial_{1}N-\partial_{2}M\\
-\partial_{1}M+\partial_{2}L
\end{bmatrix}
 +\frac1{\det \bg}
\begin{bmatrix}
 G_1 \\
G_2
\end{bmatrix},
\end{split}
\ees
where
\bes
\begin{split}
G_1=&\frac{1}{2}\big[g_{22}(-\partial _{1}g_{22})-g_{12}(\partial _{2}g_{22})\big]L
-\big[g_{22}(\partial_{2}g_{11})-g_{12}(\partial _{1}g_{22})\big]M\\
 &+\frac{1}{2}\big[g_{22}(\partial_{1}g_{11})-g_{12}(-\partial_{2}g_{11})\big]N,\\
G_2=&\frac{1}{2}\big[g_{12}(-\partial_{1}g_{22})+g_{11}(\partial_{2}g_{22})\big] L
     -\big[-g_{12}(\partial_{2}g_{11})+g_{11}(\partial_{1}g_{22})\big] M\\
&+\frac{1}{2}\big[-g_{12}(\partial_{1}g_{11})+g_{11}(-\partial_{2}g_{11})\big]N.
\end{split}
\ees
The inverse of the coefficient matrix is
$$
\begin{bmatrix}
\frac{g_{12}}{N} & \frac{g_{22}}{N} \\[1.5mm]
\frac{g_{11}}{L} & \frac{g_{12}}{L}
\end{bmatrix},
$$
which gives
\be\label{e423}
\begin{bmatrix}
\partial_{1}g_{12} \\
\partial_{2}g_{21}%
\end{bmatrix}%
=
\begin{bmatrix}
\frac{g_{12}}{N} & \frac{g_{22}}{N} \\[1.5mm]
\frac{g_{11}}{L} & \frac{g_{12}}{L}
\end{bmatrix}%
\begin{bmatrix}
 \partial_{1}N-\partial_{2}M+\frac1{\det \bg}G_1\\
-\partial_{1}M+\partial_{2}L+\frac1{\det \bg}G_2
\end{bmatrix}.
\ee

\smallskip
{\bf 4.} From \eqref{e423}, the equality of cross partials gives us the additional consistency
equation:
\be\label{e424}
\begin{split}
&\partial_{2}\left(\frac{g_{12}}{N}\Big(\partial_{1}N-\partial _{2}M+\frac1{\det \bg}G_1\Big)
+ \frac{g_{22}}{N}\Big(-\partial _{1}M+\partial _{2}L+\frac1{\det \bg}G_2\Big)\right)\\
&=\partial_{1}\left(\frac{g_{11}}{L}\Big(\partial_{1}N-\partial _{2}M+\frac1{\det \bg}G_1\Big)
+\frac{g_{12}}{L}\Big(-\partial_{1}M+\partial _{2}L+\frac1{\det \bg}G_2\Big)\right).
\end{split}
\ee
The Gauss curvature equation \eqref{e48} has the form:
\be\label{e425}
\kappa =\frac{1}{2\det \bg}\left[-\partial_{22}g_{11}+2\partial_{12}g_{12}
-\partial_{11}g_{22}\right] +l.o.t.
\ee
Substitution of \eqref{e423} into \eqref{e425} yields
the two second-order equations for $(g_{11},g_{22})$.

\smallskip
{\bf 5.} Now compute the highest order terms of the differential operators
associated with (\ref{e423}) at the origin where $g_{ij}=\delta _{ij}$:
\begin{equation}\label{e4.24-a}
\begin{split}
\begin{bmatrix}
\partial_{1}g_{12} \\[1.5mm]
\partial_{2}g_{21}%
\end{bmatrix}%
&=
\begin{bmatrix}
0 & \frac{1}{N} \\[1.5mm]
\frac{1}{L} & 0%
\end{bmatrix}%
\begin{bmatrix}
\mathbf{\ }\partial _{1}N-\partial _{2}M-%
{\frac12}%
\partial_{1}g_{22}L-\partial_{2}g_{11}M+%
{\frac12}%
\partial _{1}g_{11}N \\[1.5mm]
\mathbf{-}\partial _{1}M+\partial _{2}L+%
{\frac12}%
\partial_{2}g_{22}L-\partial _{1}g_{22}M+-%
{\frac12}%
\partial _{2}g_{11}N%
\end{bmatrix}\\[2mm]
&=%
\begin{bmatrix}
\frac{1}{N}(\mathbf{-}\partial _{1}M+\partial _{2}L+%
{\frac12}%
\partial_{2}g_{22}L-\partial_{1}g_{22}M-%
{\frac12}%
\partial_{2}g_{11}N) \\[1.5mm]
\frac{1}{L}\mathbf{\ (}\partial_{1}N-\partial_{2}M-%
{\frac12}%
\partial _{1}g_{22}L-\partial _{2}g_{11}M+%
{\frac12}%
\partial _{1}g_{11}N)%
\end{bmatrix}.
\end{split}
\end{equation}
Then \eqref{e425} becomes
\be\label{e426}
\begin{split}
\kappa =&\frac{1}{2}\left[-\partial_{22}g_{11}-%
\partial_{11}g_{22}\right]
+\frac{1}{4N}\left[
\partial_{22}g_{22}L-2\partial_{12}g_{22}M-%
\partial_{22}g_{11}N\right]\\
&+\frac{1}{4L}\left[-%
\partial_{11}g_{22}L-2\partial_{12}g_{11}M+%
\partial_{11}g_{11}N\right]+l.o.t.
\end{split}
\ee
and \eqref{e424} becomes
\begin{equation}\label{e426-b}
\begin{split}
&\frac{1}{2N}\big[
(\partial_{22}g_{22})L
-2(\partial_{12}g_{22})M-%
(\partial_{22}g_{11})N\big]\\
&=\frac{1}{2L}\big[-
(\partial_{11}g_{22})L
-2(\partial_{12}g_{11})M+%
(\partial_{11}g_{11})N\big]+l.o.t.
\end{split}
\end{equation}

In the matrix form, \eqref{e426}--\eqref{e426-b}  are written as
\be\label{e427}
\begin{bmatrix}
{\frac14}%
\frac{N}{L} & -%
{\frac14}
\\[1.5mm]
-%
{\frac12}%
\frac{N}{L} &
{\frac12}%
\end{bmatrix}%
\partial_{11}
\begin{bmatrix}
g_{11} \\[1.5mm]
g_{22}%
\end{bmatrix}%
+%
\begin{bmatrix}
-\frac{M}{2L} & -\frac{M}{2N} \\[1.5mm]
\frac{M}{L} & -\frac{M}{N}%
\end{bmatrix}%
\partial_{12}
\begin{bmatrix}
g_{11} \\[1.5mm]
g_{22}%
\end{bmatrix}%
+%
\begin{bmatrix}
-%
{\frac14}
&
{\frac14}%
\frac{L}{N} \\[1.5mm]
-%
{\frac12}
&
{\frac12}%
\frac{L}{N}%
\end{bmatrix}%
\partial_{22}
\begin{bmatrix}
g_{11} \\[1.5mm]
g_{22}%
\end{bmatrix}%
=l.o.t.
\ee
Make the change of independent variables:
$$
(x_{1}^{\prime}, x_{2}^{{\prime}})=(x_{1}+x_{2}, x_{1}-x_{2}),
$$
so that, for any $f\in C^2$,
\begin{align*}
&\frac{\partial^{2}f}{\partial x_{1}^{2}}
=\frac{\partial^{2}f}{\partial x_{1}^{\prime 2}}
+2\frac{\partial^{2}f}{\partial x_{1}^{\prime }\partial x_{2}^{\prime}}
+\frac{\partial^{2}f}{\partial x_{2}^{\prime 2}},\quad
\frac{\partial^{2}f}{\partial x_{2}^{2}}
=\frac{\partial^{2}f}{\partial x_{1}^{{\prime }2}}
-2\frac{\partial^{2}f}{\partial x_{1}^{\prime }\partial x_{2}^{\prime}}
+\frac{\partial^{2}f}{\partial x_{2}^{{\prime }2}},\\[2mm]
&\frac{\partial^{2}f}{\partial x_{1}\partial x_{2}}
=\frac{\partial^{2}f}{\partial x_{1}^{\prime 2}}
 +\frac{\partial^{2}f}{\partial x_{2}^{\prime 2}}.
\end{align*}
Then we have
\be\label{e427-a}
A\frac{\partial^{2}}{\partial x_{1}^{\prime 2}}
\begin{bmatrix}
g_{11}\\
g_{22}
\end{bmatrix}
+
A\frac{\partial^{2}}{\partial x_{2}^{\prime 2}}
\begin{bmatrix}
g_{11}\\
g_{22}
\end{bmatrix}
+\big(\frac{N}{L}+1\big)
\begin{bmatrix}
\frac{1}{2} & -\frac{1}{2} \\[1.5mm]
-1 & 1
\end{bmatrix}%
\frac{\partial^{2}}{\partial x_{1}^{\prime}\partial x_2^{\prime}}
\begin{bmatrix}
g_{11}\\
g_{22}
\end{bmatrix}
=l.o.t.
\ee
where the coefficient matrix $A$ is
\begin{align*}
A&=\begin{bmatrix}
{\frac14}%
\frac{N}{L} & -%
{\frac14}
\\[1.5mm]
-%
{\frac12}%
\frac{N}{L} &
{\frac12}%
\end{bmatrix}%
+%
\begin{bmatrix}
-\frac{M}{2L} & -\frac{M}{2N} \\[1.5mm]
\frac{M}{L} & -\frac{M}{N}%
\end{bmatrix}%
+%
\begin{bmatrix}
-%
{\frac14}
&
{\frac14}%
\frac{L}{N} \\[1.5mm]
-%
{\frac12}
&
{\frac12}%
\frac{L}{N}%
\end{bmatrix}\\[2mm]
&=
\begin{bmatrix}
{\frac14}%
\frac{N}{L}-\frac{M}{2L}-%
{\frac14}
& -%
{\frac14}%
-\frac{M}{2N}\text{ }+%
{\frac14}%
\frac{L}{N} \\[1.5mm]
-%
{\frac12}%
\frac{N}{L}+\frac{M}{L}-%
{\frac12}
&
{\frac12}%
-\frac{M}{N}+%
{\frac12}%
\frac{L}{N}%
\end{bmatrix}%
\end{align*}
with determinant equal to
$$
-{\frac12}\frac{M}{LN}\left(N+L-2M\right)
=-{\frac12}\frac{M}{LN}\left(N+L-2\sqrt{LN-\kappa}\right).
$$

If the term: $N+L-2\sqrt{LN-\kappa}=0$, then
$(N-L)^{2}=-4\kappa$.
Thus, the conditions that $\kappa >0$ and $M\neq 0$
make this coefficient matrix
non-singular.
Since  $M=-uv$ in fluid variables, we need $uv\neq 0$ and
curvature $\kappa$ determined by the Poisson equation
with $p^{2}+pq^{2}=\kappa$.

Of course, when $p>0$, $L>0$, and $N>0$, then $\frac{M}{LN}\neq 0$.
Clearly, if pressure $p$ determined by the Poisson equation is positive,
then the
condition that $\kappa >0$ is automatically satisfied.
Then line $x_1^\prime=0$ is non-characteristic for the equation in \eqref{e427-a}.

\smallskip
{\bf 6.} Since $\partial_{x_{1}^{\prime }}g_{12}={\frac12}(\partial _{x_{1}}g_{12}+\partial _{x_{2}}g_{12})$,
\eqref{e423} gives a partial differential equation for $g_{12}$ in the direction
normal to the initial data line $x_{1}^{\prime }=0$:
\begin{align}
\frac{\partial g_{12}}{\partial x_{1}^{\prime}}
=&\frac{1}{2}\big(\frac{g_{12}}{N}+\frac{g_{11}}{L}\big)\big(\partial_1N-\partial_2M+\frac{1}{\det \bg}G_1\big)\notag\\[1.5mm]
&+\frac{1}{2}\big(\frac{g_{22}}{N}+\frac{g_{12}}{L}\big)\big(-\partial_1M +\partial_2L+\frac{1}{\det \bg}G_2\big). \label{e427-b}
\end{align}

\smallskip
{\bf 7.} We now determine the initial data for system \eqref{e427-a}--\eqref{e427-b}
on line $x_{1}^{\prime}=0$ that is non-characteristic for the system.

Since
\begin{eqnarray*}
&&\partial_{x_{1}^{\prime }}g_{11}={\frac12}(\partial_{x_{1}}g_{11}+\partial _{x_{2}}g_{11}),
\quad \partial_{x_{1}^{\prime}}g_{22}={\frac12}(\partial_{x_{1}}g_{22}+\partial _{x_{2}}g_{22}),\\
&&\partial_{x_{2}^{\prime }}g_{11}={\frac12}(\partial_{x_{1}}g_{11}-\partial_{x_{2}}g_{11}), \quad
\partial_{x_{2}^{\prime }}g_{22}={\frac12}(\partial _{x_{1}}g_{22}-\partial _{x_{2}}g_{22}),
\end{eqnarray*}
the initial conditions:
$$
g_{ij}=\delta _{ij}, \quad
\partial _{x_{1}^{\prime }}g_{11}=\partial _{x_{1}^{\prime }}g_{22}=0
$$
also give us
$$
\partial_{x_{2}^{\prime}}g_{11}=\partial_{x_{2}^{\prime}}g_{22}=0.
$$
This yields
$$
\partial_{x_{1}}g_{11}=\partial_{x_{2}}g_{11}
=\partial_{x_{1}}g_{22}=\partial_{x_{2}}g_{22}=0\qquad \mbox{on $x_{1}^{\prime}=0$}.
$$
Then, from the formula:
$\partial_{x_{2}^{\prime}}g_{12}={\frac12}(\partial_{x_{1}}g_{12}-\partial_{x_{2}}g_{12})$,
\eqref{e423}, and
$$
\begin{bmatrix}
\partial_{1}g_{12} \\[1.5mm]
\partial_{2}g_{21}%
\end{bmatrix}%
=
\begin{bmatrix}
\frac{g_{12}}{N} & \frac{1}{N} \\[1.5mm]
\frac{1}{L} & \frac{g_{12}}{L}
\end{bmatrix}%
\begin{bmatrix}
\mathbf{\ }\partial_{1}N-\partial_{2}M \\[1.5mm]
\mathbf{-}\partial_{1}M+\partial_{2}L
\end{bmatrix},
$$
we obtain the ordinary differential equation \eqref{ode-1} on line $x_{1}^{\prime}=0$.
This differential equation \eqref{ode-1}
and the initial condition: $g_{12}=0$ at $x_{2}^{\prime}=0$
determine the initial data for $g_{12}$ on line $x_{1}^{\prime}=0$.
Therefore, we have specified analytic initial data
on the
non-characteristic line $x_{1}^{\prime}=0$  for system  \eqref{e427-a}--\eqref{e427-b}.

\smallskip
{\bf 8.} Finally, we see that the Cauchy-Kowalewski
theorem delivers a local analytic solution $(g_{11}, g_{12}, g_{22})$
of system \eqref{e427-a}--\eqref{e427-b}
with the analytic data
on the non-characteristic line $x_{1}^{\prime}=0$,
where we have used the analyticity of $p$ as the
solution to our Poisson equation and the analyticity of $\kappa$
from the
formula: $\kappa=p^{2}+pq^{2}$.
The positivity condition on $p$ is
irrelevant and can always be satisfied by simply adding to any $p$ solving the
Poisson equation a sufficiently large positive constant.
This completes the proof.
\endproof

Notice that the shear flow data $(u,v)=(u(x_{2}),0)$ are divergence free,
no matter how smooth or in
fact how irregular $u$ is.
In this case, $M=0$, so that Lemma \ref{L41}
does not apply to yield a metric $\bg$.
However, in this case, we have the following trivial result.

\begin{lemma}\label{L42}
If $M=0$ in any open set of $\mathbb{R}^{2}$, then there exists a developable surface
\begin{equation}\label{develop}
\y=(Ax_2, Ax_1, f(x_2))
\end{equation}
with corresponding metric $\bg^*$ given by
$$
g^*_{11}(\x)=A^2,  \quad g^*_{12}(\x)=0, \quad g^*_{22}(\x)=(f'(x_2))^2+A^2,
$$
with
$$
f'(x_2)=A\arctan\big(-A\int_0^{x_2}u^2(s) ds\big),
$$
with $A>1$ a constant.
Similarly, the corresponding developable surface with metric $\bg^*$
can be obtained for $v=v(x_1)$ as another shear flow. In particular, $g_{ij}^*(\mathbf{0})=A^2\delta_{ij}$ and $\partial_k g_{ij}^*(\mathbf{0})=0$ for $i,j,k=1,2$.

Conversely, the developable surface \eqref{develop} with $f$ and $A$ given defines a shear
flow given by \eqref{4.2-1} below and $v = 0$.
\end{lemma}

\proof $\,$ Consider the surface:
$$
y_1=Ax_2, \quad y_2=Ax_1, \quad y_3=f(x_2), \quad A=\text{const.}>1.
$$
Then
\begin{eqnarray*}
&& \partial_1 y_1=0
\quad
\partial_1 y_2=A, \quad \partial_1 y_3=0,\\
&& \partial_2 y_1=A, \quad \partial_2 y_2=0,
\quad  \partial_2 y_3=f'(x_2),
\end{eqnarray*}
so that
\begin{eqnarray*}
&& g_{11}=\partial_1 \y\cdot \partial_1 \y=A^2,\\
&& g_{12}=\partial_1 \y\cdot \partial_2 \y=0,\\
&& g_{22}=\partial_2 \y\cdot \partial_2 \y=\big(f'(x_2)\big)^2+A^2,
\end{eqnarray*}
and
$$
\mathbf{n}=\frac{\partial_1 \y\times \partial_2 \y}{|\partial_1 \y\times \partial_2 \y|}
=\frac{(f'(x_2), 0,  -A)}{\sqrt{\big(f'(x_2)\big)^2+A^2}}.
$$

Furthermore, we have
\begin{eqnarray*}
&& L=\frac{1}{\sqrt{\det{\bg}}} \partial_{11} \y\cdot \mathbf{n}=0,\\
&& M=\frac{1}{\sqrt{\det{\bg}}} \partial_{12} \y\cdot \mathbf{n}=0,\\
&& N=\frac{1}{\sqrt{\det{\bg}}} \partial_{22} \y\cdot \mathbf{n}
=\frac{1}{\sqrt{\det{\bg}}}(0,0, f''(x_2))\cdot \mathbf{n}
=-\frac{f''(x_2)}{\big(f'(x_2)\big)^2+A^2},
\end{eqnarray*}
so that
$$
\kappa=0.
$$

With $N$ prescribed by the solution of
\begin{equation}\label{4.2-1}
-\frac{f''(x_2)}{\big(f'(x_2)\big)^2+A^2}=u^2(x_2),
\qquad f'(0)=0,
\end{equation}
we have constructed the developable surface,
while the solution of \eqref{4.2-1} is determined by
$$
f'(x_2)=A\arctan\big(-A\int_0^{x_2}u^2(s)\, ds\big).
$$

Similarly, we can obtain the corresponding results for $v=v(x_1)$ as another shear flow.
\endproof

Notice that the Codazzi equations are automatically satisfied since $\Gamma_{11}^1=\Gamma_{11}^2=0$.

\begin{remark}\label{R41}
$\,$ Notice the very important fact that the shear flow initial data
given by Lemma {\rm \ref{L42}}  require only the square integrability
of the shear functions $u(x_2)$ and $v(x_1)$.
\end{remark}

\begin{remark}\label{R42}
$\,$ The assumption of analyticity is made only for
convenience. In fact, if one follows the presentation of DeTurck-Yang
{\rm \cite{DeYang}} for a similar problem where the theory of systems
of real principal type is exploited,
it is quite evident that the analyticity could be
replaced by $C^{\infty}$.
\end{remark}

\section{$\,$ Existence Theorem for the Evolving Fluids and Manifolds}\label{S5}

Now we prove our main existence theorem.

\begin{theorem}\label{T51}
For the initial data constructed in Lemma {\rm \ref{L41}}, there
exists a solution of system \eqref{e31}--\eqref{e33} and \eqref{e35},
locally in space-time,
which satisfies both the incompressible Euler equations \eqref{e25}--\eqref{e26}
and a propagating two-dimensional Riemannian manifold isometrically immersed
in $\mathbb{R}^{3}$.
\end{theorem}

\proof $\,$ For the given fluid data, by the Cauchy-Kowalewski
theorem, there exists a local space-time solution $(u,v, p)$ of
the incompressible Euler equations \eqref{e25}--\eqref{e26}
with $uv\neq 0, p>0$, and \eqref{e27}.
Define the Gauss curvature via the relation:
$$
\kappa=p^{2}+pq^{2}, \quad q^{2}=u^{2}+v^{2},
$$
so that
$LN-M^{2}=\kappa$.
Thus system \eqref{e31}--\eqref{e33} and \eqref{e35}
is satisfied if the following system can be solved:
\bes
\begin{cases}
\partial_{t}u=\Gamma_{22}^{1}L-2\Gamma_{12}^{1}M+\Gamma _{11}^{1}N,\\[1.5mm]
\partial_{t}v=\Gamma_{22}^{2}L-2\Gamma_{12}^{2}M+\Gamma _{11}^{2}N,\\[1.5mm]
\kappa \det \bg=g_{2m}(\partial_{2}\Gamma_{11}^{m}-\partial_{1}\Gamma_{12}^{m}
 +\Gamma_{11}^{n}\Gamma_{n2}^{m}-\Gamma_{12}^{n}\Gamma_{n1}^{m})
\end{cases}
\ees
for the time evolving metric $\bg$, where
\bes
\begin{cases}
u^{2}=\frac12\big[-(L-N)+\sqrt{(L-N)^{2}+4M^{2}}\big],\\[1.5mm]
v^{2}=\frac12\big[(L-N)+\sqrt{(L-N)^{2}+4M^{2}}\big],\\[1.5mm]
q^{2}=\sqrt{(L-N)^{2}+4M^{2}}.
\end{cases}
\ees
This is the same system for $\bg$
that has been solved in Lemma \ref{L41}. This completes the proof.
\endproof

\begin{corollary}\label{C51}
Let $(u,v)$ be a
solution, locally analytic in space-time, of the incompressible Euler equations
\eqref{e25}--\eqref{e26},
satisfying the initial condition $uv\neq 0$ at a point $\x_{0}$.
Then there is an evolving metric $\bg$ so that
solution $(u,v)$ of the Euler equations
also defines a time
evolving two-dimensional Riemannian manifold $({\mathcal M},\bg)$
immersed in $\mathbb{R}^{3}$.
\end{corollary}

\begin{corollary}\label{C52}
Conversely, let $({\mathcal M},\bg)$  be an isometrically embeddable
manifold in $\mathbb{R}^{3}$ with
second fundamental form $(L[\bg],M[\bg],N[\bg])$ that is
a solution of the Gauss-Codazzi equations \eqref{e21}--\eqref{e22}
{\rm (}for an as yet undetermined metric $\bg${\rm )}.
Identify the fluid variables $(u,v,p)$ by our relations{\rm :}
\be\label{5.1-a}
\begin{cases}
u^{2}=\frac{1}{2}\big[-(L-N)+\sqrt{(L-N)^{2}+4M^{2}}\big],\\[1.5mm]
v^{2}=\frac{1}{2}\big[(L-N)+\sqrt{(L-N)^{2}+4M^{2}}\big],\\[1.5mm]
p^{2}+pq^{2}=\kappa, \quad q^{2}=u^{2}+v^{2},
\end{cases}
\ee
so that $(u,v,p)=(u[\bg], v[\bg], p[\bg])$ are functionals of the unknown metric field $\bg$.
Then the equations:
\be\label{5.2-a}
\begin{cases}
\partial_{t}u = \Gamma_{22}^{1}L-2\Gamma_{12}^{1}M+\Gamma _{11}^{1}N, \\[1.5mm]
\partial_{t}v = \Gamma_{22}^{2}L-2\Gamma _{12}^{2}M +\Gamma _{11}^{2}N, \\[1.5mm]
\partial_{11}(u^{2})+2\partial_{12}(uv)+\partial_{22}(v^{2})
   =-\triangle p
\end{cases}
\ee
define a system of three equations in the three unknowns $(g_{11}, g_{12}, g_{22})$.

If this implicit system for $\bg$ has a solution, then the three functions $(u,v,p)$ satisfy
the incompressible Euler equations \eqref{e25}--\eqref{e26}.
\end{corollary}

\proof $\,$ Substitute the first two equations in \eqref{5.2-a}
into the Codazzi equations \eqref{e31}.
We see that the balance of linear momentum \eqref{e25} is satisfied,
and the third equation in \eqref{5.2-a}
implies, via the balance of
linear momentum \eqref{e25}, the incompressibility condition:
$$
\mbox{div}(u,v)=0.
$$
\endproof

\begin{remark} $\,$ The difficulty in employing Corollary {\rm \ref{C52}}
would be that the equations for $\bg$ are
only known implicitly.
Notice that equations \eqref{e33}--\eqref{e34} satisfied by $\bg$ can be written as
\be\label{e51}
\begin{cases}
 \partial_t u=\Gamma_{22}^{1}(v^{2}+p)+2\Gamma_{12}^{1}(uv)+\Gamma_{11}^{1}(u^{2}+p),\\[1.5mm]
\partial_t v=\Gamma_{22}^{2}(v^{2}+p)+2\Gamma_{12}^{2}(uv)+\Gamma_{11}^{2}(u^{2}+p),\\[1.5mm]
\partial_{1}\big(\Gamma_{22}^{1}(v^{2}+p)+2\Gamma _{12}^{1}(uv)+\Gamma_{11}^{1}(u^{2}+p)\big)\\[1.5mm]
\quad +\partial_{2}\big(\Gamma_{22}^{2}(v^{2}+p)+2\Gamma_{12}^{2}(uv)+\Gamma_{11}^{2}(u^{2}+p)\big)=0.
\end{cases}
\ee
Here $(u,v)$ are defined by \eqref{e32} in terms of $(L, M, N)$ which in turn depend
on $\bg$,  and $p$ is given by
$$
p=\frac12\left(\sqrt{q^{4}+4\kappa}-q^2\right)
$$
from \eqref{e210}.
Thus, this is a system of three equations for the three components of $\bg$.
The difficulty here is that the dependence of $(u, v)$ on $\bg$ is implicit.
It would be interesting to identify some suitable ways to check their solvability.
\end{remark}

\section{$\,$ Shear Flow Initial Data}\label{S6}

As we noted in Remark \ref{R41},
the shear flow data do not satisfy the hypothesis: $uv\neq 0$ in Lemma \ref{L41}.
Hence, this is a particularly interesting case to study.
Namely, if the fluid data of $(u,v)$ are $(u(x_{2}),0)$,
the question is whether a Riemannian
manifold $({\mathcal M},\bg)$ locally immersed in $\mathbb{R}^{3}$ can be determined.
For this data, we find immediately that $\triangle p=0$
so that $p$ is harmonic.
Furthermore, $L=p,M=0,N=u^2(x_2)+p,\kappa =p u^2(x_2)+p^{2}$,
and our constraint equations \eqref{e43}--\eqref{e44}
reduce to the equations:
\be\label{e61}
\begin{split}
&\partial_{1}p=-\Gamma_{22}^{1}p-\Gamma_{11}^{1}\big(u^2(x_2)+p\big),\\[1mm]
&\partial_{2}p=-\Gamma_{22}^{2}p-\Gamma_{11}^{2}\big(u^2(x_2)+p\big),
\end{split}
\ee
and
\be\label{e62}
\begin{split}
&pu^2(x_2)+p^{2} \notag\\
&=\frac1{2(\det \bg)^2}
\det
\begin{bmatrix}
-%
\partial_{22}g_{11}+2\partial_{12}g_{12}-%
\partial_{11}g_{22} &
\partial_{1}g_{11}
& 2\partial_{1}g_{12}-%
\partial_{2}g_{11} \\[1.5mm]
2\partial_{2}g_{12}-%
\partial_{1}g_{22} & 2g_{11} & 2g_{12} \\[1.5mm]
\partial_{2}g_{22} & 2g_{12} & 2g_{22}%
\end{bmatrix}\\[2mm]
&\quad -\frac1{2(\det \bg)^{2}}\det
\begin{bmatrix}
0 &
\partial _{2}g_{11} &
\partial _{1}g_{22} \\[1.5mm]
\partial _{2}g_{11} & 2g_{11} & 2g_{12} \\[1.5mm]
\partial _{1}g_{22} & 2g_{12} & 2g_{22}%
\end{bmatrix}.
\end{split}
\ee
As in \S \ref{S4},  \eqref{e61}--\eqref{e62}
give equation \eqref{e27} that is now
\be\label{e63}
\begin{bmatrix}
\frac{u^2(x_2)+p}{4p} & -%
{\frac14}
\\[1.5mm]
-%
\frac{u^2(x_2)+p}{2p} &
{\frac12}%
\end{bmatrix}%
\partial_{11}
\begin{bmatrix}
g_{11} \\[1.5mm]
g_{22}%
\end{bmatrix}%
+%
\begin{bmatrix}
-%
{\frac14}
&
\frac{p}{4(u^2(x_2)+p)} \\[1.5mm]
-%
{\frac12}
&
\frac{p}{2(u^2(x_2)+p)}%
\end{bmatrix}%
\partial_{22}
\begin{bmatrix}
g_{11} \\[1.5mm]
g_{22}%
\end{bmatrix}%
=l.o.t.
\ee
The symbol of the above differential operator is given by
\bes
\begin{split}
&\det \left(
\begin{bmatrix}
\frac{u^2(x_{2})+p}{4p} & -%
{\frac14}
\\[1.5mm]
-%
\frac{u^2(x_{2})+p}{2p} &
{\frac12}%
\end{bmatrix}%
\xi_{1}^{2}+%
\begin{bmatrix}
-%
{\frac14}
&
\frac{p}{4(u^2(x_{2})+p)} \\[1.5mm]
-%
{\frac12}
&
\frac{p}{2(u^2(x_{2})+p)}%
\end{bmatrix}%
\xi_{2}^{2}\right)\\[2mm]
&=\det
\begin{bmatrix}
\frac{u^2(x_{2})+p}{4p}\xi_{1}^{2}
-%
{\frac14}%
\xi_{2}^{2} & -%
{\frac14}%
\xi _{1}^{2}+%
\frac{p}{4(u^2(x_{2})+p)}\xi_{2}^{2}\\[1.5mm]
-%
\frac{u^2(x_{2})+p}{2p}\xi_{1}^{2}
-
{\frac12}%
\xi_{2}^{2} &
{\frac12}%
\xi_{1}^{2}+%
\frac{p}{2(u^2(x_{2})+p)}
\xi_{2}^{2}%
\end{bmatrix}\\[1.5mm]
&=0.
\end{split}
\ees

Thus, for the shear flow initial data, our system of partial differential equations
for $\bg$ is degenerate.
On the other hand, we know from Lemma \ref{L42}  that,
if a shear flow is given in an open neighborhood of space-time,
then
metric $\bg^*$ provides the desired map.

\medskip
\section{$\,$ The Nash-Kuiper Theorem and Existence of Wild Solutions}\label{S7}

The existence of a locally analytic Riemannian manifold with metric $\bg$
has been established in Corollary \ref{C51},
$$
g_{ij}(x_{1}, x_{2},t)
=\delta_{ij}
+h.o.t.,
$$
and hence
$(g_{ij})<(g_{ij}^*)$ in the sense of quadratic forms locally near $(x_1,x_2)=(0,0)$ for $A>1$.
Thus, we have proved that map $\y$ locally induced by metric $\bg$
is shorter than that induced by $\bg^*$  or, in language of the Nash-Kuiper
theorem \cite{CDS2012,Kuiper,Nash1954},
embedding ${\y}_\bg$ induced by $\bg$
is a
\textit{short embedding}. We state this in the following lemma:

\begin{lemma}\label{L71}
The metric, $(g_{ij}(x_{1}, x_{2},t))$,
given in Theorem {\rm \ref{T51}} and Corollary {\rm \ref{C51}}
is locally a short embedding
that satisfies
$$
(g_{ij}(x_{1},x_{2},t))< (g^*_{ij})
$$
in the sense of quadratic forms.
\end{lemma}

We then recall the Nash-Kuiper theorem as given in \cite{CDS2012}.

\begin{theorem}\label{T71}
Let $({\mathcal M}^{n}, \mathfrak{g)}$ be a smooth, compact
Riemannian manifold of dimension $n\geq 2$.
Assume that $\mathfrak{g}$ is in $C^{\infty}$.
Then

\smallskip
{\rm (i)} If $m\geq \frac{(n+2)(n+3)}{2}$,
any short embedding can be approximated by
isometric embeddings of
class $C^{\infty}$ {\rm (}cf. Nash {\rm \cite{Nash1954}} and Gromov {\rm \cite{Gromov1986}}{\rm )}{\rm ;}

\smallskip
{\rm (ii)} If $m\geq n+1,$ then any short embedding can be approximated in $C^0$ by isometric
embeddings of class $C^{1}$ {\rm (}cf. Nash {\rm \cite{Nash1954}}, Kuiper {\rm \cite{Kuiper}}{\rm )}.
\end{theorem}

We note that part (ii) of Theorem \ref{T71}
has been extended for an analytic metric $\mathfrak{g}$ with a manifold diffeomorphic to an $n$-dimensional ball
by Borisov \cite{Borisov1,Borisov2,Borisov3,Borisov4,Borisov5,Borisov6,BorisovY}
to obtain the greater regularity $C^{1,\alpha}, \alpha <\frac{1}{1+n+n^{2}}$,
and the condition of analyticity has been weakened by Conti-De Lellis-Sz\'{e}kelyhidi Jr. \cite{CDS2012}.
Without the condition that the manifold is diffeomorphic to an $n$-dimensional ball,
Conti-De Lellis-Sz\'{e}kelyhidi Jr. \cite{CDS2012} still obtain the greater
regularity $C^{1,\alpha}$ for some $\alpha >0$.

\smallskip
An immediate consequence of the results by
Nash-Kuiper \cite{Nash1954,Kuiper}, Borisov \cite{Borisov1,Borisov2,Borisov3,Borisov4,Borisov5,Borisov6,BorisovY},
and Conti-De Lellis-Sz\'{e}kelyhidi Jr. \cite{CDS2012} is the following theorem.

\begin{theorem}\label{T72}
The short embedding $\y_{\bg}$ established by Theorem
{\rm \ref{T51}} and Corollary {\rm \ref{C51}}
can be approximated in $C^0$,
by wrinkled embeddings $\y_{\rm w}$ in $C^{1,\alpha }$ for some  $\alpha \in (0,1)$, locally in space.
In elementary terms, the geometric image of our locally analytic solution
to the incompressible Euler equations, i.e.,
our surface propagating in time, can be approximated in
$C^0$ by wrinkled manifolds.
\end{theorem}

\begin{corollary}\label{c7.1} In particular, $\yy_{\rm w}(t, \cdot)$ constructed in Appendix \ref{A1} belongs
to $C^0([0,T]; C^{1,\alpha}_{\rm loc})$.
\end{corollary}

It is the time-dependent wrinkled solution of Corollary \ref{c7.1} that we use in the rest of the paper.

\begin{remark} $\,$ We reinforce the fact that, while the wrinkled embedding $\y_{\rm w}(t)$ satisfies
the time-independent relation
$\partial_{i}\y_{\rm w}\cdot \partial_{j}\y_{\rm w}=g^*_{ij}$ where $g^*_{ij}$ is independent of time,
$\y_{\rm w}$ must be time-dependent, since $\y_{\rm w}$ shadows the time-dependent short embedding $\y_{\rm g}$.

\end{remark}

\begin{theorem}\label{T74}
There are an infinite number of {\rm (}non $C^{2}${\rm )} evolving
manifolds for metric $\mathfrak{g}_{ij}=g^*_{ij}$,
all with the same
initial wrinkled data at $t=0$.
Furthermore, for all these solutions, the energy remains constant{\rm :}
$$
E(t)=\int\limits_{\Omega}(\partial_{i}\y_{\rm w} \cdot \partial_i\y_{\rm w}){\rm d}\x =const.
$$
where the manifold is defined by
the $C^{1,\alpha}$--solution
of $\partial_{i}\y_{\rm w}\cdot \partial_{j}\y_{\rm w}=g^*_{ij}$
which trivially satisfies $\partial_{i}\y_{\rm w}\cdot \partial_{i}\y_{\rm w}=g^*_{ii}$
so that $\partial_{i}\y_{\rm w}$ is in $L^{\infty}$ in time.
\end{theorem}

\proof $\,$ On any sequence of time intervals,
just switch back and forth
from any of the infinite choices of wrinkled solutions.

Specifically, let $[0,T]$ be a time interval in which there exists a unique
analytic solution to the Cauchy
initial value problem for the Euler equations.
By Theorem \ref{T72}, this solution is the pre-image of a short embedding $\y_{\bg}$ and
can be approximated by a wrinkled embedding $\y_{\rm w}$.

Divide the interval $[0,T]$ into sub-intervals
$[T_0,T_1], [T_1, T_2], \dots, [T_{n-1}, T_n]$ with $T_0=0$ and $T_n=T$.
Take any sequence $\varepsilon_k>0, k=1,\cdots, n$.
Then, by Theorem \ref{T72}, we have a sequence of
wrinkled solutions $\{\y^k_{\rm w}\}$ such that
$$
\|\y_{\bg}-\y^k_{\rm w}\|_{C}<\varepsilon_k
\qquad \mbox{for $T_{k-1}\le t\le T_k$, $k=1, \dots, n$}.
$$
Now define the wrinkled solution $\y^*_{\rm w}$ on $[0,T]$ by
$$
\y^*_{\rm w}:=\y_{\rm w}^k
\qquad \mbox{for $T_{k-1}\le t< T_k$, $k=1, \dots, n$}.
$$

Fix $\y_{\rm w}^1$ as the Cauchy data,
but allow $\varepsilon_k$ and $T_k$ to vary for $k=2, \dots, n$.
In this way, we have produced an infinite number of wrinkled solutions satisfying the same Cauchy data.
Then the solutions are in $C^{1,\alpha}$ in space and continuous in time on
each of the sub-intervals.

To compute the
energy, use
$$
E(t)=\int\limits_{\Omega }(\partial_{i}\y_{\rm w}\cdot \partial_{i}\y_{\rm w})\,{\rm d}\x
=\int\limits_{\Omega }g^*_{ii}\,{\rm d}\x =const.
$$
\endproof

In the results given in Theorems \ref{T71}--\ref{T74},
by construction,
embedding $\y_{\rm \bg}$ is
the geometric image of a smooth solution of the Euler
equations \eqref{e25}--\eqref{e26}.
On the other hand, the Nash-Kuiper $C^{1,\alpha}$ embeddings have
not been shown to be the image of solutions in any sense of the Euler
equations \eqref{e25}--\eqref{e26}.
We now provide this link.

Let $\Omega\Subset \mathcal{M}$ be a compact domain, where $\mathcal{M}$ is a regular surface with a family
of Riemannian metric $\{\bg(t)\}_{t\in[0,T]}$.
Consider the difference $\vv(t):=\y_{\rm w}(t)-\y_{\rm \bg}(t) \in C^{1,\alpha}(\Omega; \R^3)$ for each $t \in [0,T]$
and continuous in $t$.
As computed by G\"{u}nther (see \cite{gunther}), one obtains
\begin{equation}\label{A}
\partial_i (\vv\cdot\partial_j\y_{\rm \bg})+\partial_j (\vv\cdot \partial_i \y_{\rm \bg})
- 2\vv \cdot \big(\Gamma^k_{ij} \partial_k \y_{\rm \bg} + H_{ij}\mathbf{n}\big) = h_{ij},
\end{equation}
where $\Gamma^k_{ij}$ is the Christoffel symbol of $(\mathcal{M}, \bg)$, $H \in {\rm Sym}^2(T^*\mathcal{M})$
is the second fundamental form associated with the smooth isometric embeddings
$\y_{\rm \bg}(t)$, and $\n$ is the outward unit normal.
Moreover, we suppress subscript $t$ in the equations which are {\em kinematic} ({\it i.e.}, pointwise in $t$).
We introduce the symmetric quadratic form  $\bh=(h_{ij})\in C^0( {\rm Sym}^2(T^*\mathcal{M}))$ by
\begin{equation}\label{h}
h_{ij} := g^*_{ij}-{g}_{ij}-\partial_i \vv \cdot \partial_j \vv.
\end{equation}

Now we {\em project} the difference $\vv$ along the direction $\partial_i \y_{\rm \bg}$:  Define
\begin{equation}
\bar{v}(t)^i := \vv(t) \cdot \partial_i \y_{\rm \bg}(t) \in C^{1, \alpha}(\Omega; \R)
\qquad \text{ for each } t \in [0,T], \, i \in \{1,2\}.
\end{equation}
This is valid since $\y_{\rm \bg}$ is smooth. Noting that, for any vector field $\phi=(\phi^j)$,
one has
\begin{equation*}
\nabla_i \phi^j = \partial_i \phi^j - \Gamma^k_{ij} \phi^k,
\end{equation*}
with $\nabla$ being the covariant derivative on $(\mathcal{M}, \bg)$.
Thus, \eqref{A} can be recast into
\begin{equation}\label{B}
\nabla_i \bar{v}^j + \nabla_j \bar{v}^i - 2\vv \cdot H_{ij}\n = h_{ij}.
\end{equation}

Multiplying \eqref{B} by $A^{kij}$,
we obtain a system of two first-order  scalar PDEs:
\begin{equation}\label{key equation}
\bar{A}^1 \partial_1 \bar{\vv} + \bar{A}^2 \partial_2 \bar{\vv} + B\bar{\vv} = \bar{\bh},
\end{equation}
where
\begin{equation}\label{regularity}
\begin{cases}
\bar{\vv}=(\bar{v}^1, \bar{v}^2)^\top \in C^{1,\alpha}(\Omega; \R^2),\\
\bar{\bh}=(A^{1ij}h_{ij}, A^{2ij}h_{ij})^\top \in C^0(\Omega; \R),\\
B= (B_{ij})= (-A^{jlm} \Gamma^i_{lm})\in C (\mathfrak{gl}(2;\R)),
\end{cases}
\end{equation}
and the $2\times 2$ matrices
\begin{equation}
\bar{A}^k := (A^{ijk})_{1\leq i,j \leq 2} \qquad \text{ for } k \in \{1,2\},
\end{equation}
satisfy the symmetry property: $A^{kij}=A^{kji}=A^{ikj}$,
Notice that, due to the symmetry,
one has
$$
\bar{A}^1=\begin{bmatrix}
A^{111}& A^{112}\\
A^{112}&A^{122}
\end{bmatrix},\quad
\bar{A}^2=\begin{bmatrix}
A^{112}&A^{122}\\
A^{122}&A^{222}
\end{bmatrix},
$$
matrices $\bar{A}^k$ have only four independent components
$(A^{111}, A^{112}, A^{122}, A^{222})$.

Now we consider   $(U_1,U_2,P)$ satisfying  the incompressible Euler equations:
\begin{equation}\label{B3}
\begin{split}
&\partial_tU_1+\partial_1\left(U_1^2+P\right)+\partial_2\left(U_1U_2\right)=0,\\
&\partial_tU_2+\partial_1\left(U_1U_2\right)+\partial_2\left(U_2^2+P\right)=0,\\
&\partial_1U_1+\partial_2U_2=0.
\end{split}
\end{equation}
On the other hand, \eqref{key equation} may be written as
\begin{equation}\label{B2}
\partial_1\left(\bar{A}^1\begin{bmatrix} \bar{v}_1 \\ \bar{v}_2\end{bmatrix}\right)
+\partial_2\left(\bar{A}^2\begin{bmatrix} \bar{v}_1 \\ \bar{v}_2\end{bmatrix}\right)
-\partial_1\bar{A}^1\begin{bmatrix} \bar{v}_1 \\ \bar{v}_2\end{bmatrix}
-\partial_2\bar{A}^2\begin{bmatrix} \bar{v}_1 \\ \bar{v}_2\end{bmatrix}
=\begin{bmatrix}  {h}_1 \\  {h}_2\end{bmatrix},
\end{equation}
where $\mathbf{h}=\begin{bmatrix}  {h}_1 \\  {h}_2\end{bmatrix}=\bar{\mathbf{h}}-B\bar{\vv}$.
Thus, the identification:
\begin{equation}\label{e710}
\begin{split}
&U_1^2+P=A^{111}\bar{v}_1+A^{112}\bar{v}_2=:r_1,\\
&U_1U_2=A^{112}\bar{v}_1+A^{122}\bar{v}_2=:r_2,\\
&U_2^2+P=A^{122}\bar{v}_1+A^{222}\bar{v}_2=:r_3,
\end{split}
\end{equation}
will force the second two terms of the first two equations in \eqref{B3} to agree with
$$
\partial_1\left(\bar{A}^1
\begin{bmatrix} \bar{v}_1 \\ \bar{v}_2
\end{bmatrix}\right)
+\partial_2\left(\bar{A}^2
\begin{bmatrix} \bar{v}_1 \\ \bar{v}_2
\end{bmatrix}\right).
$$
Furthermore, it is a simple matter to express $(U_1, U_2, P)$ in terms of $(r_1, r_2, r_3)$
as
\begin{equation}\label{e711}
U_1=(r_1-P) ^\frac12, \quad
U_2=(r_3-P)^\frac12,
\end{equation}
where $P$ satisfies
\begin{equation}\label{e712}
P=\frac12\left((r_1+r_3)-\left((r_1-r_3)^2+4r_2^2\right)^\frac12\right),
\end{equation}
so that $r_1-P>0$, $r_2-P>0$.

Thus, for \eqref{key equation} to agree with the incompressible Euler equations, we see the system:
\begin{align}
&\partial_t\begin{bmatrix}U_1\\ U_2\end{bmatrix}
 +\partial_1\bar{A}^1\begin{bmatrix}\bar{v}_1\\ \bar{v}_2\end{bmatrix}
 +\partial_2\bar{A}^2\begin{bmatrix}\bar{v}_1\\ \bar{v}_2\end{bmatrix}
 +\bh=0,   \label{e713}\\[1.5mm]
&\partial_1U_1+\partial_2U_2=0, \label{e714}
\end{align}
must be satisfied, where $(U_1, U_2)$ depend smoothly on $(A^{111}, A^{112}, A^{122}, A^{222})$,
and $\bh$ depends linearly on $A^{kij}$.
This is a under-determined system, since there are three equations in the four unknowns
$(A^{111}, A^{112}, A^{122}, A^{222})$.
A non-trivial solution of this system thus maps the Nash-Kuiper solution to a weak solution
of the incompressible Euler equations.

Let us rewrite \eqref{e714} in terms of a stream function $\Psi(\x,t)$:
\begin{align}
&U_2+\partial_1\Psi=0, \label{e715}\\
&U_1-\partial_2\Psi=0, \label{e716}
\end{align}
where, from \eqref{e711},
$$
U_i=\hat{U}_i(A^{111}, A^{112}, A^{122}, A^{222}; v_1, v_2), \qquad i=1,2.
$$
Define $\ff:=(A^{111}, A^{112}, A^{122}, A^{222},\Psi)$. Then \eqref{e715} is of the form:
\begin{equation}\label{e717}
\Phi_1(\x,t,\ff,\partial_1\ff)=0,
\end{equation}
while \eqref{e716} is of the form:
\begin{equation}\label{e718}
\Phi_2(\x,t,\ff,\partial_1\ff,\partial_2\ff)=0,
\end{equation}
where $\Phi_1$ and $\Phi_2$ are scalar valued functions.
Furthermore, \eqref{e713} is of the form
\begin{equation}\label{e719}
\Phi_3(\x,t,\ff,\partial_t\ff,\partial_1\ff,\partial_2\ff)=0,
\end{equation}
where $\Phi_3$ takes values in $\R^2$.

Thus we have written our system \eqref{e713}--\eqref{e714} for $\ff$
in Gromov's triangular form \eqref{e717}--\eqref{e719}.
In fact, let us quote Gromov's theorem \cite[p. 198]{Gromov1986} verbatim as follows:
{\it Let $\partial_i, i=1, \dots, k$, be continuous linearly independent vector fields on $V$.
For example, $V=\R^n$ and $\partial_i=\partial/\partial u_i$, $i=1,\dots,n$.
Let $\Phi_i$ be smooth vector valued functions such that $\Phi_i$ takes values in $\R^{s_i}$ and $\Phi_i$
has entries $\vv, \ff, \partial_1 \ff, \dots, \partial_i\ff$, where $\ff$ is the unknown map $V\to\R^q$.
In other words, $\Phi_i: V\times\R^{q(i+1)}\to\R^{s_i}$.
Consider the following {\rm (}triangular{\rm )} systems of $s=\sum_{i=1}^ks_i$ PDEs{\rm :}
\begin{equation}\label{e720}
\begin{split}
&\Phi_1(\vv,\ff,\partial_1\ff)=0,\\
&\Phi_2(\vv,\ff,\partial_1\ff,\partial_2\ff)=0,\\
&\cdots \cdots \cdots\\
&\Phi_k(\vv,\ff,\partial_1\ff,\dots,\partial_k\ff)=0.
\end{split}
\end{equation}
{\bf Local Solvability.} If $s_i\le q-1$ for all $i=1,\dots,k$, and if the functions $\Phi_i, i=1,\dots, k$, are generic,
then system \eqref{e720} admits a local $C^1$--solution $\ff: U\to\R^q$, for some open subset $U\subset V$.
Moreover, the $C^1$-solutions $U\to\R^q$ are $C^0$--dense in some open subset in the space of $C^0$-maps $U\to\R^q$.}

Notice that Gromov's theorem requires ``smoothness" and ``genericity" of $(\Phi_1, \Phi_2, \Phi_3)$.
The ``genericity"  appears to be a requirement of nonlinearity on $(\Phi_1, \Phi_2, \Phi_3)$
which is satisfied because of the nonlinear relation \eqref{e711}.
However, ``smoothness" is not obvious, since $\Phi_i$ is at most $C^{0,\alpha}$ in space
due to the occurrence of derivatives of $V$.
Since $\y_{\rm w}$ is continuous in $t$, but not necessarily smooth, ``smoothness"   in $t$
is also an issue.
Nevertheless, if Gromov's theorem is still valid in our case, we would obtain an infinite number of
solutions $\ff$. However, from \eqref{e715}--\eqref{e716},
$\ff$ in $C^1$ would  yield $(U_1, U_2)$ at best continuous
in $\x=(x_1, x_2)$, and hence a weak but not strong solution of the incompressible Euler equations.
Hence, it seems that, for the moment, we have formal but not yet rigorous map from
our Nash-Kuiper solution to the Euler equations.

\section{$\,$ The Compressible Euler Equations}\label{S8}

The arguments we have used for the incompressible case also carry over to
the compressible Euler equations.

Then the equations for the balance of linear momentum are
\be\label{e81}
\begin{cases}
\partial_{1}(\rho u^{2}+p)+\partial_{2}(\rho uv)=-\partial_{t}(\rho u),\\[1.5mm]
\partial_{1}(\rho uv)+\partial_{2}(\rho v^{2}+p)=-\partial_{t}(\rho v),
\end{cases}
\ee
The equation for the conservation of mass is given by the equation:
\be\label{e82}
\partial_{1}(\rho u)+\partial_{2}(\rho v)=-\partial_{t}\rho.
\ee
Set
\be\label{e83}
L=\rho v^{2}+p,\quad M=-\rho uv,\quad N=\rho u^{2}+p,
\ee
and the Gauss equation becomes
$$
(\rho v^{2}+p)(\rho u^{2}+p)-(\rho uv)^{2}=\kappa,
$$
{\it i.e.},
\be\label{e84}
p^{2}+p\rho q^{2}=\kappa,
\ee
where $q^{2}=u^{2}+v^{2}$.
For the compressible case, we take
\be\label{e85}
p=\frac{\rho^{\gamma}}{\gamma}, \qquad \gamma\ge 1.
\ee
Substitute \eqref{e85} into \eqref{e84} to obtain
\be\label{e86}
\Big(\frac{\rho}{\gamma}\Big)^{2\gamma}+ q^2\frac{\rho^{\gamma +1}}{\gamma}=\kappa.
\ee

For simplicity, let us first take the isothermal case $\gamma =1$ so that
\be\label{e87}
\rho^{2}(1+ q^2)=\kappa,
\ee
that is,
\be\label{e88}
\rho =\sqrt{\kappa/(1+q^{2})}.
\ee
From \eqref{e88}, $\rho$ as the density
is determined explicitly as a function of $\kappa$ and $q^{2}$.
On the other hand, $\rho$ must satisfy the mass balance equation
\eqref{e82} so that
\be\label{e89}
\partial_{1}(u\sqrt{\kappa/(1+q^{2})})+\partial_{2}(v\sqrt{\kappa/(1+q^{2})})
= -\partial_{t} (\sqrt{\kappa/(1+q^{2})}).
\ee
We see that, for given $(u, v)$,
\eqref{e89} is a scalar conservation law for the Gauss
curvature $\kappa$.
This leads us to the following elementary lemma:

\begin{lemma}\label{L81}
For a given smooth solution $(\rho, u, v)$ of the
compressible Euler equations \eqref{e81}--\eqref{e82},
the initial value problem for the Gauss
curvature $\kappa$ given by \eqref{e89} with initial data
$\kappa=\rho^{2}(1+q^{2})>0$ at $t=0$
has a global  smooth solution $\kappa$ in space-time,
which satisfies \eqref{e87}.
\end{lemma}

\proof $\,$  Notice that equation \eqref{e89} is a linear transport
equation of conservation form for $\frac{1}{\sqrt{\kappa}}$.
Then the result holds simply by using the standard local
existence-uniqueness theorem for the conservation law
and the direct computation for $J=\sqrt{\frac{\kappa}{1+q^2}}-\rho$:
\begin{align*}
\partial_{t}J&=-(\partial_{1}\rho u+\partial_{2}\rho v)+
\partial_{1}(u\sqrt{\kappa/(1+q^{2})})
+\partial_{2}(v\sqrt{\kappa/(1+q^{2})})\\[1mm]
&=-\partial _{1}(Ju)-\partial_{2}(Jv).
\end{align*}
This equation has the unique
solution $J=0.$
\endproof

The general case $\gamma \geq 1$ follows from noting that
the left-hand side of \eqref{e86} is a monotone increasing function
of $\rho$ so that, for any $\kappa>0$,
there exists a function $\rho(\kappa, q)$ that solves \eqref{e86}.
The regularity of this solution is seen by differentiation
of \eqref{e86} with respect
to $(q, \kappa)$.
We can thus state the general version of Lemma \ref{L81} as

\begin{lemma}\label{L82}
$\,$For a given smooth solution $(\rho, u, v)$ of $\,$the
compressible Euler
equations  \eqref{e81}--\eqref{e82},
the initial value problem for the Gauss
curvature  for \eqref{e89} with initial data \eqref{e86}
at $t=0$ has a global solution in space-time
satisfying \eqref{e86}.
\end{lemma}

Again, we can also write the fluid variables in terms of the geometric
variables $(L, M, N)$.
To do this, simply substitute \eqref{e85}
with $\rho$ given by the solution to \eqref{e86}
into $L=\rho v^{2}+p, M=-\rho uv$, and $N=\rho u^{2}+p$ to
find
$$
L-N=\rho \big(v^{2}-u^{2}\big).
$$
Write $v=-\frac{M}{\rho u}$ to see
$L-N=\rho \big(\frac{M}{\rho u}\big)^{2}-\rho u^{2}$, which yields
$$
(L-N)\rho u^{2}=M^{2}-\rho ^{2}u^{4},\quad
\rho ^{2}u^{4}+(L-N)\rho u^{2}-M^{2}=0.
$$
Then
$$
u^{2}=\frac{1}{2\rho}\big[-(L-N)\pm \sqrt{(L-N)^{2}+4M^{2}}\big],
$$
and, with $L-N=v^{2}-\big(\frac{M}{v}\big)^{2}$,
$$
v^{2}=\frac{1}{2\rho}\big[(L-N)\pm \sqrt{(L-N)^{2}+4M^{2}}\big].
$$
This tells us to choose sign $``+ "$ in the above formulas so that
\be\label{e810}
\begin{cases}
u^{2}=\frac{1}{2\rho}\big[-(L-N)+\sqrt{(L-N)^{2}+4M^{2}}\big],\\[1.5mm]
v^{2}=\frac{1}{2\rho}\big[(L-N)+\sqrt{(L-N)^{2}+4M^{2}}\big],
\end{cases}
\ee
which are the desired relations.

Now the compressible analogs of Theorem \ref{T51} and Corollaries \ref{C51}--\ref{C52}
follow by the same arguments we have employed in \S \ref{S5}.

\begin{theorem} \label{T81}
For locally analytic initial data $(\rho, u, v)$
satisfying \eqref{e86} and $\rho uv\neq 0$,
there exists a local solution in space-time
of system \eqref{e41b}--\eqref{e45b} which satisfies both the compressible Euler
equations \eqref{e81}--\eqref{e82} and a
propagating two-dimensional Riemannian manifold  isometrically immersed in $\mathbb{R}^{3}$.
\end{theorem}

\begin{corollary}\label{C81}
Let $(\rho, u, v)$ be a local analytic
solution in space-time of the compressible Euler equations \eqref{e81}--\eqref{e82}
satisfying the initial condition{\rm :}
$\rho uv\neq 0$ at a point $x_{0}$.
Then there is an evolving metric $\bg$ so
that solution $(\rho, u, v)$ of the compressible Euler equations
\eqref{e81}--\eqref{e82}
also
defines a time evolving two-dimensional Riemannian manifold $({\mathcal M}, \bg)$ immersed
in $\mathbb{R}^{3}$ with the Gauss curvature
$\kappa$
determined \eqref{e86}.
\end{corollary}

Theorems \ref{T72}--\ref{T74} remain unchanged, except for the fact that
metric $\bg^*$ for the wrinkled manifold now corresponds to the
special incompressible solution of the compressible Euler equations \eqref{e81}--\eqref{e82}
with
$\rho=const.\, $
In particular, Theorem \ref{T72} now gives the smooth solutions of the
geometric image of the compressible Euler equations \eqref{e81}--\eqref{e82}
being approximated by
the wrinkled solutions that correspond
weak
shear solutions of
the \emph{incompressible} Euler
equations.
The construction of the wrinkled solutions in this case
is the same as done in \S \ref{S7}
and produces the solutions of the incompressible Euler equations \eqref{e25}--\eqref{e26}.
A wrinkled solution of the compressible Euler equations \eqref{e81}--\eqref{e82}
is impossible, since it would correspond to a vacuum via (\ref{e84}).

\begin{remark} $\,$ We note that the results given in \S 7 and this section on
the existence of {\it wild}
weak shear solutions have been given in terms of
the Cartesian coordinates $\x=(x_{1}, x_{2})$.
 The choice of local
Cartesian coordinates
is only a convenience, and the Euler
equations written in polar coordinates would suffice.
\end{remark}

\section{$\,$ Isometric Embedding Problem and General Continuum Mechanics: Elastodynamics as an Example}\label{S10}

Motivated by our results for fluid dynamics,
we now consider solutions in two-dimensional general continuum mechanics.
Denote by $\T=(T_{11}, T_{12}, T_{22})$ the (symmetric) Cauchy stress tensor,
and assume that fields $(u,v)$, $\T$,
and $\rho$ are consistent with some specific constitutive equation
for a body and satisfy the balances of mass and linear momentum
(satisfaction of the balance of angular momentum
is automatic).

The equations for the balance of linear momentum in the spatial representation are
\be\label{e101}
\begin{cases}
\partial_{1}(\rho u^{2}-T_{11})+\partial_{2}(\rho uv-T_{12})
  =-\partial_{t}(\rho u),\\[1mm]
\partial_{1}(\rho uv-T_{12})+\partial_{2}(\rho v^{2}-T_{22})
 =-\partial_{t}(\rho v),
\end{cases}
\ee
and the balance of mass is
\begin{equation}\label{9.2a}
\partial_{1}(\rho u) +\partial_{2}(\rho v) = -\partial_{t}\rho,
\end{equation}
or
\begin{equation}\label{9.2b}
\rho = \rho_0 \, (\det \F)^{-1},
\end{equation}
where $\rho_0$ is the density of the body in the reference configuration,
and $\F$ is the deformation gradient of the current configuration
with respect to this reference.

\subsection{$\,$ Mapping a general continuum mechanics problem
to the non-degenerate isometric embedding problem}

Denote the geometric dependent variables as
\be\label{LMN}
L=\rho v^{2}-T_{22},\quad N=\rho u^{2}- T_{11},\quad
M= -\rho uv + T_{12},
\ee
so that the Gauss equation becomes
\begin{equation}\label{10.5-a}
(\rho v^{2}-T_{22})(\rho u^{2}-T_{11})-(\rho uv-T_{12})^{2}= \kappa.
\end{equation}

Under the assumption that
\begin{equation}\label{10.6-a}
\det \T+ 2\rho uvT_{12}-\rho v^{2}T_{11}-\rho u^{2}T_{22} > 0,
\end{equation}
at least on some initial time interval,
such a solution of continuum mechanics corresponds to a positive Gauss curvature
for the corresponding isometric embedding problem.
Using $(L, M, N)$ as time-dependent data in (\ref{e31})--(\ref{e32})
and expressing these three equations in terms of the metric components $(g_{11}, g_{12}, g_{22})$
and their derivatives by using (\ref{e23})--(\ref{e24}),
we have a system of three partial differential equations for the components of the metric.
Using the required analog of Corollary \ref{C81}, then we have

\begin{corollary}
Let $(u,v, T_{11}, T_{12}, T_{22}, \rho)$ be a local analytic  solution in space-time
of system \eqref{e101}--\eqref{9.2a}
at a point satisfying the condition{\rm :}
$$
T_{12}-\rho uv  \neq 0\qquad \mbox{at a point $\x^0$}.
$$
Then there is an evolving metric $\bg$ so that the continuum mechanical solution
defines a time evolving two-dimensional Riemannian manifold $({\mathcal M},\bg)$
isometrically immersed in $\mathbb{R}^3$.
\end{corollary}

\subsection{$\,$ Image of degenerate isometric embedding problem in continuum mechanics}
We note that the ability to solve for the evolving metric $\bg$ in the fluid case
relied on the initial data for the off-diagonal term $M$ in the second fundamental
form being non-vanishing.
In the cases of incompressible and compressible fluids,
this term was simply $M=-uv$.
However, for general continuum mechanics,
$$
M = T_{12}-\rho u  v,
$$
that is, the expression for $M$ has an additional contribution.
Theorems \ref{T72}--\ref{T74} remain unchanged, except for the fact that metric $\bg^*$
for the wrinkled manifold now corresponds
to special \emph{steady} smooth solutions of the equations of two-dimensional general continuum mechanics.
We now identify these solutions.

We want to define solutions to the mechanical equations \eqref{e101}--\eqref{9.2a}
in continuum mechanics
from a smooth, degenerate isometric embedding problem for
which $\partial_i \y \cdot \partial_j \y = g^*_{ij}, i,j=1,2$, are satisfied.
Then, making the association analogous to that used in Lemma 4.2,
we have
\be\label{10.8-a}
\begin{cases}
\rho uv - T_{12} = 0,\\
\rho u^2 - T_{11} = Z(x_2,t),\\
\rho v^2 - T_{22} = 0,
\end{cases}
\ee
where function $Z$ is defined from the embedding, and
\be
\begin{cases}
\partial_1(\rho uv - T_{12}) + \partial_2 (\rho v^2 - T_{22}) = 0, \\
\partial_1( \rho u^2 - T_{11}) + \partial_2( \rho uv - T_{12}) = 0
\end{cases}
\ee
are satisfied by the mechanical fields that are being defined.
These would form a consistent set of fields satisfying the balances
of linear momentum and mass if the following constraints
\be\nonumber
\begin{cases}
 \partial_t (\rho u) = \partial_t (\rho v) = 0,\\
 \partial_1 (\rho u) + \partial_2 (\rho v) = -\partial_t \, \rho
\end{cases}
\ee
are satisfied, {\it i.e.},
these geometric solutions define continuum mechanical solutions with steady momenta.
This is easily done by noting that the conditions imply $\partial_{tt} \rho = 0$ with the solution
\begin{equation}\label{rho_flat}
\rho (x_1, x_2, t) = \rho_1(x_1,x_2)\, t + \rho_2 (x_1,x_2),
\end{equation}
where $\rho_1$ and $\rho_2$ are arbitrary time-independent functions
of the spatial variables, which are required to be so chosen that $\rho>0$.
We then define $(u, v)$ by integrating the pointwise ordinary differential equations:
\be\nonumber
\begin{cases}
\rho u_t = - \rho_1 u,\\
\rho v_t = - \rho_1 v.
\end{cases}
\ee
With $(\rho, u, v, Z)$ in hand, we define the stress components
from \eqref{10.8-a} to obtain
a class of mechanical solutions in general continuum mechanics
from the smooth degenerate isometric embedding problem.
It is to be noted that, for such continuum mechanical solutions to be realizable
for a specific material, the mechanical fields as defined have to be shown
to be consistent with a constitutive equation for the stress for that material.

Alternatively, by seeking solutions to the system
\be\label{steady}
\begin{cases}
\partial_1(\rho uv - T_{12}) + \partial_2 (\rho v^2 - T_{22}) = 0, \\
\partial_1( \rho u^2 - T_{11}) + \partial_2( \rho uv - T_{12}) = 0,\\
\partial_1 (\rho u) + \partial_2(\rho v) = -\partial_t \,\rho,
\end{cases}
\ee
we can define solutions to the balance of mass and the \emph{steady equations of balance of linear momentum}.

\begin{remark} $\,$ The solutions to system \eqref{steady} do not necessarily constitute exact, steady solutions
of the balance of linear momentum, {\it i.e.}, $\partial_t (\rho u)$ and $\partial_t (\rho v)$ may not evaluate to $0$
from such motions, much in the spirit of quasi-static evolutions in solid mechanics and Stokes flow
in fluid mechanics.
As is understood, such solutions are typically interpreted in an asymptotic sense when the velocities
are assumed to ``equilibrate'' on a much faster time-scale than the evolution of driving boundary conditions
or forcing, {\it i.e.}, $t$ is assumed to be the ``slow'' time scale and the right-hand sides of the first two
equations actually have the terms: $-\epsilon \, \partial_t (\rho u)$ and $-\epsilon \, \partial_t (\rho v)$,
with $0 < \epsilon \ll 1$, so that we deal with the singular perturbation with a first-order approximation.
\end{remark}

\smallskip
We now
display the image of the isometric embedding problem
in the smooth degenerate case in nonlinear elastodynamics for a Neo-Hookean material.

We choose the constitutive equation for the Cauchy stress of the compressible material as
\[
\T = \rho \F\F^\top =: \rho \B,
\]
where $\F$ is the deformation gradient from a fixed reference configuration of the body.
The generic point on the fixed reference is denoted by $\X = (X_1, X_2)$, the motion as $\x(\X,t)$,
and $F_{ij} = \frac{\partial x_i}{\partial X_j}$.
For the sake of simplicity, we ignore a multiplicative scalar function of the invariants of $\B$
that would ensure that the stress at the reference configuration vanishes, even in this simplest
frame-indifferent elastic constitutive assumption.

\emph{We now seek special solutions to the steady equations \eqref{steady}.
After obtaining any such solution consistent with the posed steady problem,
we will further check which of them also constitute {\rm (}steady{\rm )} solutions to the equations of balance
of linear momentum}.

\smallskip
Thus, it needs to be demonstrated that the conditions:
\be\label{NH_gov_eq}
\begin{cases}
 uv - B_{12} = 0,\\
 v^2 - B_{22} = 0,\\
 u^2 - B_{11} = \rho^{-1} Z(x_2,t)
\end{cases}
\ee
are satisfied, along with balance of mass in the form
\[
\rho = (\det \F)^{-1} \rho_0(\X),
\]
where $\rho_0$ is the mass density distribution in the reference configuration.
We assume the mass density distribution on the reference to be a constant function with value $\rho_0 > 0$, and always require that $\det \F >0$.

The components of $\B$ are given as
\[
B_{11} = F_{11}^2 + F_{12}^2, \quad B_{12} = B_{21} = F_{11}F_{21} + F_{12}F_{22},
\quad B_{22}=F_{21}^2 + F_{22}^2.
\]
Thus, the equations required to be satisfied by a motion consistent with these constraints
would be
\be\label{elastic_syst}
\begin{cases}
\frac{\partial x_1}{\partial t}\frac{\partial x_2}{\partial t}
 = \frac{\partial x_1}{\partial X_1}\frac{\partial x_2}{\partial X_1}
 + \frac{\partial x_1}{\partial X_2}\frac{\partial x_2}{\partial X_2},\\[1.5mm]
\Big(\frac{\partial x_2}{\partial t}\Big)^2
  = \Big(\frac{\partial x_2}{\partial X_1}\Big)^2
   + \Big(\frac{\partial x_2}{\partial X_2}\Big)^2,\\[1.5mm]
\Big(\frac{\partial x_1}{\partial t}\Big)^2
  = \Big(\frac{\partial x_1}{\partial X_1}\Big)^2
   + \Big(\frac{\partial x_1}{\partial X_2}\Big)^2 + \left(\rho_0^{-1} \det \F \right) Z (x_2(X,t),t).
\end{cases}
\ee
We define $\nabla x_i := (\partial_{X_1} x_i, \partial_{X_2} x_i)$,
and
\begin{align*}
J (\nabla x_1, \nabla x_2):=& \det \F = \partial_{X_1} x_1 \partial_{X_2} x_2 - \partial_{X_2} x_1 \partial_{X_1} x_2\\[1mm]
=& \sqrt{(\nabla x_1 \cdot \nabla x_2)^2 - |\nabla x_1|^2 |\nabla x_2|^2}
\end{align*}
for subsequent use.

If solutions exist to the above system,
then they also satisfy the following conditions:
the satisfaction of the first equation of (\ref{elastic_syst}) by a solution
satisfying the second and third equations is equivalent to
\[
\rho_0^{-1} Z |\nabla x_2|^2 + J(\nabla x_1, \nabla x_2) = 0,
\]
which implies that $|\nabla x_2| > 0$, or else $ J = 0$ which is not acceptable by hypothesis.

Thus, the solutions of the steady, Neo-Hookean image of the degenerate isometric embedding problem must satisfy the following system for functions $(x_1, x_2)$:
\be\label{es2}
\begin{cases}
\Big(\frac{\partial x_2}{\partial t}\Big)^2 = |\nabla x_2|^2,\\[1mm]
 \Big(\frac{\partial x_1}{\partial t}\Big)^2
 = |\nabla x_1|^2 - \Big(\frac{J\left(\nabla x_1, \nabla x_2 \right)}{|\nabla x_2|}\Big)^2,\\
 \rho_0^{-1} Z(x_2, t) |\nabla x_2|^2 = - J(\nabla x_1, \nabla x_2),
\end{cases}
\ee
with the caveat that the solutions to \eqref{es2} satisfy
$$
\frac{\partial x_1}{\partial t} \frac{\partial x_2}{\partial t} =  \nabla x_1 \cdot \nabla x_2.
$$

When restricting attention to the class of solutions to the whole system that are scale-invariant
in $(x_1, x_2)$, {\it i.e.}, those solutions that remain solutions if the dependent and independent
variables are scaled by the same constant,
we note that the first two equations in (\ref{es2}) remain invariant under the rescaling of
the independent
and dependent variables by $\lambda$.
Similarly, $Z(\lambda x_2, \lambda t)$ is independent of $\lambda$.
\emph{We now assume that the given function $Z(x_2,t) = - \rho_0$}.

\begin{example}\label{elast_example}
$\,$ Consider a shearing motion of the form:
\be\label{ex_ansatz}
\begin{split}
& x_1 = X_1 + f(X_2,t),\\
& x_2 = X_2 \pm t.
\end{split}
\ee
Then $F_{11} = 1$, $F_{12} = \partial_{X_2} f$, $F_{21} = 0$, and $F_{22} = 1$.
Thus, $J = 1$, $|\nabla x_1|^2 = 1 + (\partial_{X_2} f)^2$, and $|\nabla x_2|^2 = 1$.
Then the first equation of (\ref{es2}) is identically satisfied, the second requires
\[
\Big(\frac{\partial f}{\partial t}\Big)^2 = \Big(\frac{\partial f}{\partial X_2} \Big)^2,
\]
and the third is identically satisfied by our choice of $Z = -\rho_0$.

Thus, traveling waves of the form:
\[
f(X_2, t) = w(X_2 \pm t)
\]
define solutions to system \eqref{es2}
for each sign of $\pm$ in the second equation of (\ref{ex_ansatz}).

For a spatially uniform density field on the reference configuration,
it is easy to check that the equations of balance of linear momentum for the Neo-Hookean constitutive assumption
we are considering (the first Piola-Kirchhoff stress is given by $\rho_0 \F$) reduces to the linear,
second order system:
\[
\frac{\partial^2 x_i}{\partial t^2}
= \frac{\partial^2 x_i}{\partial X_1^2} + \frac{\partial^2 x_i}{\partial X_2^2}.
\]
Thus,  for the assumed ansatz (\ref{ex_ansatz}), the image of the degenerate isometric embedding problem
has produced special, \emph{exact} non-steady, solutions to the balance of linear momentum.
Clearly, these are not steady solutions on the reference configuration in general.
Thus, it is an interesting question to check whether any of these
are \emph{exact} steady solutions on the current configuration.

\begin{enumerate}
\item[(i)] The solutions generated in Example \ref{elast_example} correspond to
$f( {{X}_{2}},t)=w( {{X}_{2}}\pm t )$.
Let us consider each in turn.
To check steadiness, we consider an arbitrarily fixed point
$\x^0=(x_{1}^{0},x_{2}^{0})$ in space.
We have to show that the velocities of material points that occupy it
at different times is the same,
since the density is constant, $\rho={{\rho }_{0}}$, everywhere for this example.
Let the image in the reference configuration of point $\x^{0}$
at time $t$ be
\[
(X_{1}^{0}(\x^{0},t ), X_{2}^{0}(\x^{0},t)).
\]
For $f({{X}_{2}},t)=w({{X}_{2}}+t)$,
\begin{equation}
 x_{1}^{0}=X_{1}^{0}+w(X_{2}^{0}+t), \quad  x_{2}^{0}=X_{2}^{0}+t,
\end{equation}
which implies
\begin{equation}\label{10.14a}
X_{1}^{0}=x_{1}^{0}-w( x_{2}^{0}).
\end{equation}

Now
\begin{equation}
u({{X}_{1}},{{X}_{2}},t)=\frac{\partial f}{\partial t}({{X}_{2}},t )=w'( {{X}_{2}}+t), \quad
v( {{X}_{1}},{{X}_{2}} )=1.
\end{equation}
Then
\begin{equation}
\begin{split}
&u( X_{1}^{0},X_{2}^{0},t )=w'( X_{2}^{0}+t )=w'( x_{2}^{0}-t+t )=w'( x_{2}^{0}),\\
&v( X_{1}^{0},X_{2}^{0},t)=1.
\end{split}
\end{equation}
Thus, we indeed have a steady spatial velocity and momentum field.
The way to  understand it physically is the following:
Let $w(X_2)$ be specified.
The points on line ${{x}_{2}}=x_{2}^{0}$ have material points sitting on them at time $t$
with referential $X_2$--coordinate given by $x_{2}^{0}-t$.
However, for such referential points, the horizontal velocity $w'$ corresponds
to
$$
X_2+t=x_{2}^{0}-t+t=x_{2}^{0}.
$$
Hence, this is a very curious situation where the picture is completely steady on the current configuration;
however,
there is unsteady wave propagation with Piola-Kirchhoff shear
stress waves propagating on the reference configuration.

\smallskip
\item[(ii)]
Now we consider the other solution $f({{X}_{2}},t)=w( {{X}_{2}}-t)$. Then
	\begin{equation}\nonumber
	\begin{split}
  & X_{2}^{0}=x_{2}^{0}-t, \\
 & X_{1}^{0}=x_{1}^{0}-w( X_{2}^{0}-t)=x_{1}^{0}-w( x_{2}^{0}-2t),
\end{split}
\end{equation}
so that
\begin{equation}\nonumber
   u( X_{1}^{0},X_{2}^{0} )=-w'( x_{2}^{0}-2t), \quad
    v( X_{1}^{0},X_{2}^{0})=1.
\end{equation}
Then we have unsteadiness in the current configuration as well as the reference.
Note that this case corresponds to a situation that does not satisfy the caveat on
solutions
of \eqref{es2} to correspond to solutions of \eqref{elastic_syst}.

\smallskip
\item[(iii)] However, if we now choose the ansatz for the motion to be
	\begin{equation*}
  {{x}_{1}}={{X}_{1}}+f( {{X}_{2}},t), \quad {{x}_{2}}={{X}_{2}}-t,
\end{equation*}
then it can be checked that again $f({{X}_{2}},t )=w( {{X}_{2}}\pm t)$
are both solutions to system \eqref{es2},
but now only solution $w( {{X}_{2}}-t)$ works (for obvious reasons by following the previous argument):
\begin{equation*}
\begin{split}
X_{2}^{0}&=x_{2}^{0}+t, \\
X_{1}^{0}&=x_{1}^{0}-f(X_{2}^{0},t) = x_{1}^{0}-w( X_{2}^{0}-t)\\
&=x_{1}^{0}-w( x_{2}^{0}+t-t )=x_{1}^{0}-w( x_{2}^{0}).
\end{split}
\end{equation*}
Then again
\begin{equation*}
u(X_{1}^{0},X_{2}^{0})=w'( X_{2}^{0}-t)=w'( x_{2}^{0}), \quad
v( X_{1}^{0},X_{2}^{0})=v( {{X}_{1}},{{X}_{2}})=-1.
\end{equation*}
Thus, we have a steady field on the current configuration, while, on the reference,
stress/velocity/position waves are moving from the bottom to the top.
\end{enumerate}

\smallskip
While simple, it is important to appreciate that the generated exact solution
in this extremely simple example can produce smooth analogs (with continuous deformation gradient)
of the traveling wave profile shown in Fig. \ref{f1}
\begin{figure}
\centering
\includegraphics[width=5.0in, height=3.0in]{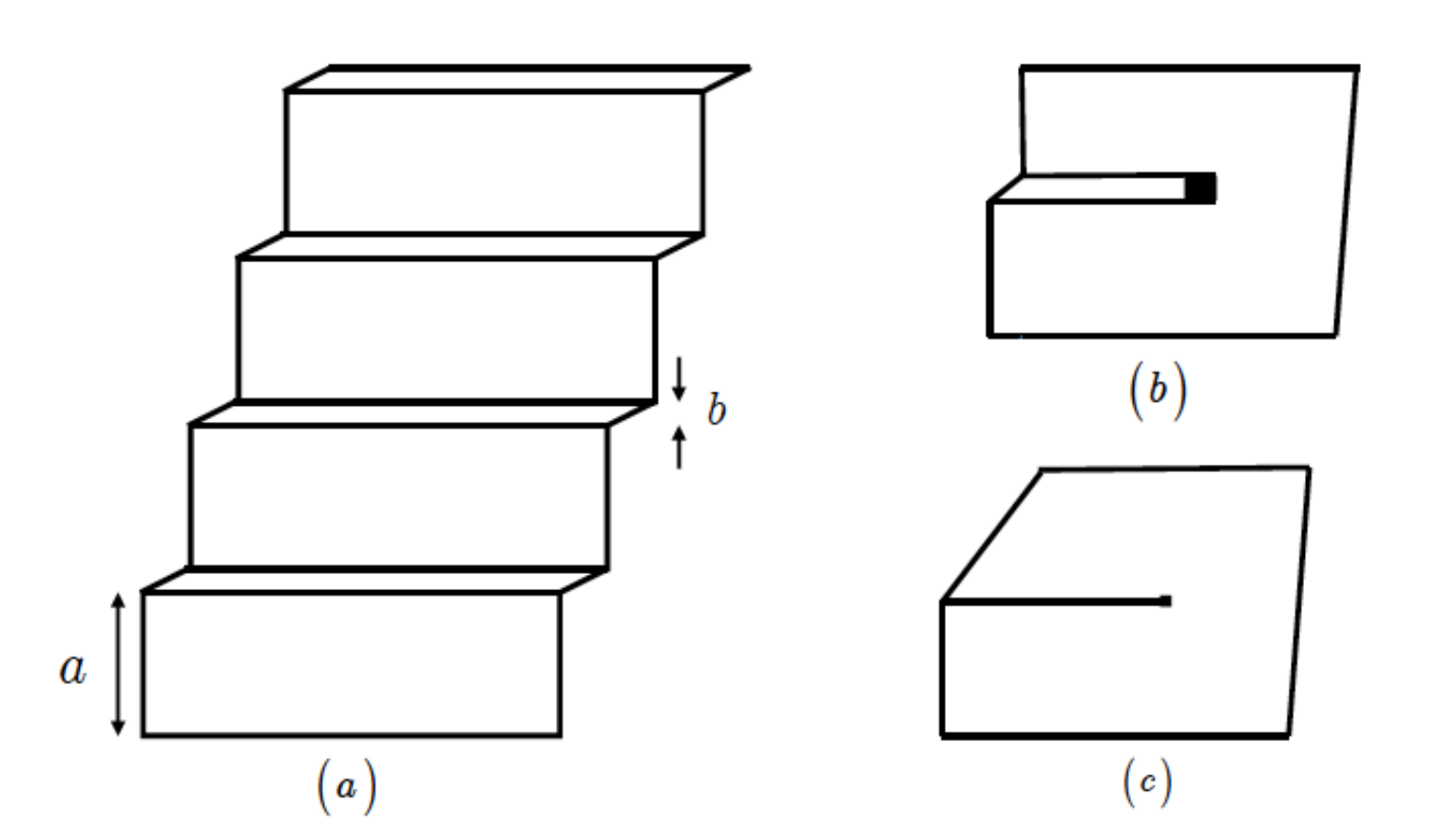}
\caption{(a) Shear bands and phase boundaries;
(b) Dislocation with core, the terminating line (through the plane of paper)
of a shear band;
(c) Generalized disclination with core, the terminating line of a phase boundary.}
\label{f1}
\end{figure}
with arbitrarily small positive values of $a$ and $b$.
When $a$ is comparable to $b$,
this is a sequence of phase boundaries separating domains;
when $b \ll a$, then the regions spanned within width $b$ are slip zones or shear bands
separating undeformed blocks.
The greater freedom that elastodynamics provides over elastostatics in producing microstructure
can now be appreciated.
While  such deformations in elastostatics can only be produced with a multi-well energy
function \cite{Ball_James, Abeyaratne_Knowles},
such microstructures (for $a,b \ll 1$) in elastodynamics occur
in this case of the simplified Neo-Hookean
material (with the energy density $tr(\mathbf{C})$),
and persist for all times. Indeed, multiple low-energy states corresponding to large total strains
and static microstructures\footnote{Note that our example is unequivocally dynamic owing to the requirement that $|\nabla x_2|^2 > 0$.}
are observed facts, and their representation consists an important physical ingredient of solid mechanics\footnote{Indeed, in the traditional definition
of phase transformations in solids, {\it cf.} \cite{Abeyaratne_Knowles}, the microstructures being discussed here would not be considered as domains
separated by phase boundaries, since the strain states belong to the same well.}.
Our intent here is to simply point out the greater freedom of displaying microstructure without any length scale
in elastodynamics.
By considering the case for ever smaller $a,b > 0$, it is also now perhaps easy to intuitively see why the wrinkled
deformations of the isometric embedding problem may have a limiting status in elastodynamics,
as the parameters $a$ and/or $b$ tend to zero.
Furthermore, it is an interesting question whether our elastodynamic model (\ref{es2}) corresponding to the degenerate
isometric embedding also displays solutions that represent static and dynamic terminating lines of phase
boundaries (generalized disclinations) and shear bands (dislocations),
ingredients of nature that may help to further understand physically the nature of the wrinkled embeddings
from geometry.
Finally, it is important to recognize and accept that microstructure in solids is not about infinite refinement
and, as is well-known, such scale-free deformations, no matter how exciting from the analytical point of view,
are notoriously difficult to deal with practically in modeling endeavors (following nature's guide, in some sense).
Thus, our considerations seem to point to the need for models that can represent both static, quasi-static
({\it i.e.}, evolving microstructure in the absence of material inertia), as well as fully dynamic microstructures
with in-built length scales and accounting for the dissipation due to microstructure evolution.
\end{example}

\subsection{$\,$ Mapping between the general continuum mechanics problem and the isometric embedding problem}
Suppose that we have a smooth, time-dependent two-dimensional manifold immersed in $\mathbb{R}^3$.
Then $(g_{11}, g_{12}, g_{22})$ and $(L,M,N)$ can be defined from the manifold
as in \S \ref{S2}.
Consider the initial data to be available on $(\rho, u, v)$ denoted by $(\rho_0, u_0, v_0)$ pointwise.
Then define the following quantities $(m_u, m_v)$ through the ordinary differential equations:
\be\label{momenta}
\begin{cases}
\partial_{t}\, m_u = \Gamma_{22}^{1}L-2\Gamma_{12}^{1}M+\Gamma_{11}^{1}N,\\
\partial_{t}\, m_v = \Gamma_{22}^{2}L-2\Gamma_{12}^{2}M+\Gamma_{11}^{2}N,
\end{cases}
\ee
with initial condition on $m_u$ specified as $\rho_0 u_0$ and on $m_v$ as $\rho_0 v_0$.
With fields $(m_u, m_v)$ available, solve
\begin{equation}\label{bal_mass}
-\partial_t \, \rho = \partial_1 m_u + \partial_2 m_v
\end{equation}
with initial condition on $\rho$ to be $\rho_0$.
Then, using this $\rho$ field along with $\Gamma$ and $(L, M, N)$ as time-dependent data,
solve the pointwise ordinary differential equations for fields $(u, v)$ given by (\ref{momenta})
with $(m_u, m_v)$ replaced by $(\rho u, \rho v)$, with the initial condition on $(u, v)$ being $(u_0, v_0)$, respectively.
Clearly, we have
\[
(m_u, m_v) = (\rho u, \rho v).
\]
However, then (\ref{bal_mass}) implies that the constructed fields $(\rho, u,v)$ satisfy the balance of mass.
Finally, defining the stress components from (\ref{LMN}) by using $(L, M, N)$ and the constructed fields $(\rho, u, v)$ as the data,
and noting (\ref{momenta}) and (\ref{e31}), we find that every smooth time evolving
two-dimensional manifold immersed in $\mathbb{R}^3$ defines a solution of the balance laws of two-dimensional
general continuum mechanics.
It requires the imposition of further constraints to obtain the motions within this class that are consistent
with the constitutive response of any specific material.

For a given constitutive relation, the Cauchy stress is determined from  the velocity field $(u,v)$ and density $\rho$.
The case for (incompressible and compressible) inviscid fluids has been exposited,
and the case of elasticity requires the deformation gradients to be calculated from the deformation field
inferred from the velocity field $(u,v)$ via time integration.
Hence, since $(\rho,u,v)$ have  been determined by our argument above,
the ability to satisfy the constitutive relations requires that metric $\bg$ be
consistent with these extra constraints.

Let $\F$ be the deformation gradient field of the mechanical body from a fixed stress-free elastic
reference configuration.
Let $\LL_v$ denote the velocity gradient field.
For simplicity,  assume a local constitutive equation for the stress of the form:
\[
\T = \hat \T (\LL_v,\rho,\F,\bg).
\]
We also use the notation: $\kappa = \kappa(\bg,\nabla \bg, \nabla^2 \bg)$
and $\Gamma= \Gamma(\bg,\nabla \bg)$, for the known functions for  the Gauss curvature and the Christoffel symbols.

\begin{remark} $\,$ A time-dependent set of mechanical fields $(\rho, u, v, \F, \T)$
and a time-dependent metric field $\bg$ related
by the constitutive assumption $\hat \T$  are consistent in the sense that the balance of linear momentum and
balance of mass are satisfied, and the metric is isometrically embedded in $\mathbb{R}^3$ at each instant of time,
provided that the following system of constrained partial differential equations of evolution are satisfied:
\begin{equation}\label{constrained_system}
\begin{cases}
\partial_{1}(\rho u) +\partial_{2}(\rho v) = -\partial_{t}\rho,\\[1mm]
\partial_k F_{ij} v_k - (\LL_v)_{ik}F_{kj} = -\partial_t F_{ij},\\[1mm]
\partial_{1}\big(\rho u^{2}-\hat T_{11}(\LL_v,\rho, \F, \bg)\big)
+\partial_{2}\big(\rho uv- \hat T_{12}(\LL_v,\rho, \F, \bg)\big)
 =  -\partial_{t}(\rho u),\\[1mm]
\partial_{1}\big(\rho uv - \hat T_{12}(\LL_v,\rho, \F, \bg)\big)
 +\partial_{2}\big(\rho v^{2} - \hat T_{22}(\LL_v,\rho, \F, \bg)\big)=-\partial _{t}(\rho v),
\end{cases}
\end{equation}
and
\begin{equation}\label{constrained_system-a}
\begin{cases}
\partial_{1} N -\partial_{2}M
 = -\Gamma_{22}^{1}(\bg,\nabla \bg) L +2\Gamma_{12}^{1} (\bg,\nabla \bg) M -\Gamma_{11}^{1} (\bg,\nabla \bg)N,\\[1mm]
\partial_{1} M -\partial_{2}L
  = \Gamma_{22}^{2}(\bg,\nabla \bg) L -2\Gamma_{12}^{2}(\bg,\nabla \bg) M +\Gamma_{11}^{2}(\bg,\nabla \bg) N,\\[1mm]
 LN - M^2 = \kappa(\bg,\nabla \bg, \nabla^2 \bg),\\[1mm]
L \mathbf{e}^1 \otimes \mathbf{e}^1
+ M \left(\mathbf{e}^1 \otimes \mathbf{e}^2 + \mathbf{e}^2 \otimes \mathbf{e}^1 \right)
+ N \mathbf{e}^2 \otimes \mathbf{e}^2\\
$\qquad$ = \mathbf{f} \left( \rho \mathbf{u} \otimes \mathbf{u}
- \mathbf{T}\left(\mathbf{L}_v, \rho, \mathbf{F}, \mathbf{g} \right) \right),
\end{cases}
\end{equation}
where $\mathbf{f}$ is a tensor-valued function of its tensorial argument,
and $\left(\mathbf{e}^1, \mathbf{e}^2\right)$ is the dual basis corresponding
to the natural basis on the surface given
by $\left(\partial_1\mathbf{y},\partial_2 \mathbf{y} \right)$.
The last tensorial equation consists of three independent equations.

The equations in \eqref{constrained_system-a}
are to be considered as the constraints
that determine the family of pairs of first and second fundamental forms of
embedded manifolds consistent with a given mechanical state at any given time.
{\it Conversely}, the evolution of the mechanical fields following
the first four equations \eqref{constrained_system}
must be constrained to the {\it manifold} defined by
\eqref{constrained_system-a}
in the state-space of spatial $(\rho, u, v, \F, \bg)$ fields.

One may consider eliminating variables $(L, M, N)$ from the equations in \eqref{constrained_system-a}
to obtain three constraint equations for the three components of the metric field.

Abstractly, one may think of eliminating all of the equations in \eqref{constrained_system-a}
and replacing $\bg$ in the mechanical set of the first four equations \eqref{constrained_system}
as a spatially non-local term in the mechanical fields representing a solution of \eqref{constrained_system-a}.

Notice that this remark applies as well to the steady problem of continuum mechanics
(where the right-hand sides of the third and fourth equations of \eqref{constrained_system}
are assumed to be $0$) and the time-dependence of boundary conditions
(or body forces that have been assumed to vanish here for simplicity)
drive the evolution of the mechanical problem.

We note that the earlier sections of this paper dealing with the equations
of incompressible and compressible fluid dynamics,
and Neo-Hookean elastodynamics are special cases of the above system where the constitutive equation
for the Cauchy stress are independent of metric $\bg$.
\end{remark}

\section{$\,$ Concluding Remarks}\label{S9}

We close by discussing some broad implications and possible extensions of the presented work.

\subsection{$\,$ Admissibility of weak solutions}
Since the dynamics of the wrinkled solutions shadowing the
incompressible or compressible Euler equations are completely reversible
Nash-Kuiper solutions corresponding to the metric for the developable surface, the usual
irreversible entropy admissibility criteria become useless as selection
criteria. This tells us that the only hope of selecting a unique physically
meaningful solution would be chosen without recourse to time evolution, {\it e.g.},
artificial viscosity or energy minimization.
For example, an energy
minimization which penalizes second derivatives of $\y$ would prefer the affine
initial data over data with folds, so that the {\it wild} initial data would be ruled out.
Similarly,
Sz\'{e}kelyhidi Jr. \cite{SzL2011}, Bardos-Titi-Wiedermann \cite{BTitiW},
Bardos-Titi \cite{BTiti}, and Bardos-Lopes Filho-Niu-Nussenzveig Lopes-Titi \cite{BLFN}
have shown that the viscosity criterion also eliminates the {\it wild} solutions for the
Euler equations.

To be more precise about the role of viscosity, we recall the incompressible
Navier-Stokes equations:
\be\label{e91}
\begin{cases}
\partial_{1}(u^{2}+p)+\partial_{2}(uv)=-\partial_{t}u+\frac{1}{\text{Re}}\Delta u,\\[1mm]
\partial_{1}(uv)+\partial _{2}(v^{2}+p)=-\partial_{t}v+\frac{1}{\text{Re}}\Delta v,
\end{cases}
\ee
where we have taken density $\rho =1$ and $\text{Re}$ denotes the
Reynolds number.

The condition of incompressibility and the Poisson equation for the pressure
remain as
\begin{align}
&\partial_{1}u+\partial_{2}v=0,\label{e92}\\[1mm]
&\partial_{11}^{2}(u^{2})+2\partial_{12}(uv)+\partial_{22}(v^{2})
=-\Delta p. \label{e93}
\end{align}
If we review all the previous arguments made for the inviscid Euler
equations leading up to and including \S \ref{S8},
we see all the conclusions
we have made regarding the Euler equations hold true for the Navier-Stokes
equations modulo one crucial point.
For the Euler equations,
metric $\bg^*$ provides the map from the steady shear, $p=0$,
a solution of the Euler equations to the Gauss-Codazzi equations.
However, for a shear solution of the Navier-Stokes equations,
the right-hand sides of \eqref{e91}
must vanish.
Therefore, instead of the fluid pre-image being a steady shear,
the fluid pre-image must satisfy the diffusion equation:
\[
-\partial_{t}u+\frac{1}{\text{Re}}\Delta u=0,\quad v=0,\quad p=0,
\]
or
\[
-\partial _{t}v+\frac{1}{\text{Re}}\Delta v=0,\quad u=0,\quad p=0.
\]
Thus, the pre-image is smooth so that, by \eqref{e28},
we have a smooth second fundamental form and a smooth embedding $\y$.
However, this contradicts the non-$C^{2}$ property of
our Nash-Kuiper wrinkled solution.
Thus, the only possibility is that, for the Navier-Stokes equations,
there is
no fluid pre-image of the Nash-Kuiper wrinkled solutions.
One may be
tempted to discount the physical relevance of the Nash-Kuiper theorem
outlined here and in the recent work in
\cite{BDIS, BDS1, CDK, CDS2012, DS2009, DS2010, DS2012, DS2013, DS2014, DS2015,SzL2011,Wied}.
On the other hand,
as noted in Chen-Glimm \cite{CG2012}
and others ({\it cf}. \cite{R-K,Wal}),
viscous fluid turbulence arises with the imposition of an external force,
which overcomes the viscous dissipation.
Hence, it may be that the addition of an external force would
indeed bring us back to the Nash-Kuiper-Gromov turbulence scenario.
Furthermore, the dynamics of defect microstructure like dislocations,
phase and grain boundaries,
triple junctions, and point defects in crystalline solids furnish a compelling
physical argument for accepting/developing physically rigorous and practically
computable models that account for the representation of microstructure,
necessarily then not of infinite refinement and with a modicum of uniqueness
in the predicted evolution of their fields.
In analogy with a crumpled piece of paper that does not produce infinitely fine
terminated folds and neither unfolds itself back to its original flat state,
perhaps the correct notion of admissibility is to move to augmented physical models
of continuum defect dynamics involving extra kinematics representing the microscopic,
and hence smoothed, dynamics of discontinuity surfaces, their terminating lines,
and point singularities of the fields of the original macroscopic model
(like nonlinear elasticity and Navier-Stokes),
while accounting for the energetics of these defects and the dissipation produced
owing to their motion. Such partial differential equation-based augmentations
of nonlinear elasticity theory have begun to emerge, {\it e.g.},
Acharya-Fressengeas \cite{ach_fress}, along with their interesting predictions
of soliton-like dynamical behavior of nonsingular
defects (Zhang et al. \cite{zhangetal}) and their collective behavior.

\subsection{$\,$ A plausible minimal model for internally stressed elastic materials}
The constrained evolution system \eqref{constrained_system}--\eqref{constrained_system-a},
written on the deforming configuration of an elastic body being tracked in Lagrangian fashion,
appears to pose an interesting model for internally stressed elastic bodies
whose constitutive response in terms of the deformation gradient
and a metric representing a stress-free state is known:
\begin{equation}\label{elast_constrained_system}
\begin{cases}
\rho = \rho_0 (\det \F)^{-1},\\[1mm]
\partial_{1}\hat T_{11}(\F, \bg)+\partial_{2}\hat T_{12}(\F, \bg)=\rho\, d_t u,\\[1mm]
\partial_{1}\hat T_{12}(\F, \bg)+\partial_{2}\hat T_{22}(\F, \bg)= \rho \,d_t v,
\end{cases}
\end{equation}
and
\begin{equation}\label{elast_constrained_system-a}
\begin{cases}
 \partial_{1} N -\partial_{2}M
 = -\Gamma_{22}^{1}(\bg,\nabla \bg) L +2\Gamma_{12}^{1} (\bg,\nabla \bg) M -\Gamma_{11}^{1} (\bg,\nabla \bg)N,\\[1mm]
\partial_{1}M -\partial_{2}L = \Gamma _{22}^{2}(\bg,\nabla \bg) L -2\Gamma_{12}^{2}(\bg,\nabla \bg) M +\Gamma_{11}^{2}(\bg,\nabla \bg) N,\\[1mm]
 LN - M^2 = \kappa(\bg,\nabla \bg, \nabla^2 \bg),\\[1mm]
 L \mathbf{e}^1 \otimes \mathbf{e}^1
 + M \left(\mathbf{e}^1 \otimes \mathbf{e}^2 + \mathbf{e}^2 \otimes \mathbf{e}^1 \right)
 + N \mathbf{e}^2 \otimes \mathbf{e}^2
 = \mathbf{\hat{f}} \left(  \mathbf{T}\left(\mathbf{F}, \mathbf{g} \right) \right),
\end{cases}
\end{equation}
where $\partial_i, i=1,2$, represent the spatial derivatives on the current configuration,
$d_t$ represents the material time derivative operator,
$\mathbf{\hat{f}}$ is a tensor-valued function of its tensorial argument,
and $\left(\mathbf{e}^1, \mathbf{e}^2\right)$ is the dual basis corresponding to the natural
basis on the surface given by $\left(\partial_1\mathbf{y},\partial_2 \mathbf{y} \right)$.
The last tensorial equation consists of three independent equations.

First of all, we note that, on dimensional grounds,
the identification of the second fundamental form with mechanical objects with physical dimensions
of stress implies from the Codazzi
equations that the metric is physically
dimensionless and the Christoffel symbols have dimensions of reciprocal lengths.
Thus, the metric may be interpreted as describing strain.
Of course, this identification also implies that a material parameter with physical units
of ($\mbox{stress}\cdot\mbox{length})^2$ is required to make the Gauss curvature equation
dimensionally consistent, while introducing a length-scale into the traditional elastic problem
of internal stress.

Next, we consider the equations in (\ref{elast_constrained_system-a}) for given $\F$.
Considering, for the moment, the situation when the constitutive equation is independent of $\bg$,
this becomes a question of determining the metric, given the second fundamental form,
such that an embedding exists in $\mathbb{R}^3$,
which is the opposite of the isometric embedding problem that may be interpreted as the question
of determining the second fundamental form, given a metric $\bg$.
At any rate, as some of our results show, this problem has a solution in many instances.
It is also perhaps reasonable to expect that the situation does not change drastically
even when the constitutive equation depends on $\bg$,
and is not degenerate in the sense that, for almost all $\F$,
there are many solutions for $\bg$ (this would necessitate evolution equations for $\bg$).
Based on this premise, the requirement of an embedding of $\bg$ in $\mathbb{R}^3$ assumes
a physical status replacing a separate constitutive equation for the  evolution of $\bg$
that would be required for the mechanical problem otherwise.
Moreover, much like in systems displaying relaxation oscillations,
the geometry and stability of states on the constraint manifold can,
on occasion, lead to interesting dynamical behavior of the mechanical fields.
The above model appears to have the possibility of being considered as a minimal model
for the statics and dynamics of soft and biological
materials ({\it e.g.}, Efrati et al
\cite{kupferman},
Jin et al
\cite{suo}, and Ambrosi et al \cite{goriely}) depending only on the knowledge of the elastic
response of the material and the material constant needed to define the Gauss curvature equation.

\begin{appendices}
\section{Time-Continuity of the Wrinkled Solutions}\label{A1}

In this appendix,  let us prove the  following result:

Let $(u,v,p)$ be a solution to the two-dimensional Euler equations in the neighbourhood
of a point $\x^\star \in \mathcal{M}$, locally analytic in space, such that
\begin{equation}\label{cauchy-k condition}
u(\x^\star,0) v(\x^\star,0) \neq 0 \qquad \text{ at } t = 0.
\end{equation}
Moreover, let $\yy_\bg(t,\cdot): (\mathcal{M}, \bg) \hookrightarrow \R^3$ be the corresponding short immersion
of an analytic surface for each $t\in [0,T]$, with $(\mathcal{M}, \bg)$ being
the Riemannian manifold given by Theorem 5.1 and Corollary 5.1.

Denote by
\begin{equation*}
\{\Phi_t\}_{t\in [0,T]}: X \equiv C^\infty_{\rm loc}(\mathcal{M}; \R^3)
\longrightarrow Y \equiv C^{1,\alpha}_{\rm loc}(\mathcal{M}; \R^3)
\end{equation*}
the collection of maps sending the short immersion $\yy_\bg(t,\cdot)$ to the {\it wild} isometric
immersion $\yy_w(t,\cdot)$, indexed by time $t$,
constructed following Conti-De Lellis-Sz\'{e}kelyhidi Jr. \cite{CDS2012}.
It is defined as follows:
\begin{align}\label{arrows}
\Phi_t: \quad X\equiv C^\infty_{\rm loc}(\mathcal{M}; \R^3)
&\xrightarrow{ \Phi_1 } C^1_{\rm loc}(\mathcal{M}; \R^3)
\xrightarrow{ \Phi_2 } C^2_{\rm loc}(\mathcal{M}; \R^3)\nonumber\\
&\xrightarrow{ \Phi_3 } C^{1,\alpha}_{\rm loc}(\mathcal{M}; \R^3) \equiv Y,
\end{align}
where $\uu:=\Phi_1(\yy_\bg(t, \cdot))$ is the $C^1_{\rm loc}$ isometric immersion
constructed by Nash-Kuiper in \cite{Kuiper, Nash1954},
$\vv:=\Phi_2(\uu)$ is the $C^2_{\rm loc}$ map ``close to being isometric'',
obtained by mollifying $\uu$, and $\yy_{\rm w}:=\Phi_3(\vv)$ is the $C^{1,\alpha}_{\rm loc}$ isometric immersion
constructed by Conti-De Lellis-Sz\'{e}kelyhidi Jr. \cite{CDS2012}.

\begin{theorem}
For each $\yy_\bg\in (X, C^1_{\rm loc})$, we have
\begin{equation*}
\Phi(\yy_\bg( \cdot))\in C^0\left([0,T], C^{1}_{\rm loc}\right).
\end{equation*}
Therefore£¬ if $u(\x,t_k)v(\x,t_k)\neq 0$ locally in space and
$(u,v,p)(t_k)\to (u,v,p)(t)$ in $C^1_{\rm{loc}}$ as $t_k\to t$,
then
$\yy_{\rm w}(t_k)\to \yy_{\rm w}(t)$ in $C^1_{\rm{loc}}$ as $t_k\to t$.
\end{theorem}

\begin{remark} $\,$ First, analytic solutions are obtained via the Cauchy-Kowalewski theorem,
which entails condition \eqref{cauchy-k condition}.
Second, the wrinkled solutions $\yy_{\rm w}$ are constructed from the short map $\yy_\bg$
by adding ``Nash wrinkles'' or ``corrugations'',
whose first derivatives are of only H\"{o}lder regularity at best.
The current upper bound, obtained by Borisov in \cite{Borisov6} and by Conti-De Lellis-Sz\'{e}kelyhidi Jr.
in \cite{CDS2012}, is $\alpha < \frac{1}{7}$£¬
where $\frac{1}{7}=\frac{1}{1+2J_n}$ for $n=2$, with $J_n = \frac{n(n+1)}{2}$
known as the {\em Janet dimension}.
\end{remark}

\proof $\,$ We divide the proof into five steps.

\smallskip
{\bf 1.} {\it Reduction to a geometric problem}.
First, we reduce the problem to showing the $C^1$--continuity of
wrinkled solutions $\yy_{\rm w}$ with respect to the short maps $\yy_\bg$.

Recall from \S 2 that the Gauss curvature is defined from the fluid variables as
\begin{equation*}
\kappa:=p^2 + p(u^2+v^2).
\end{equation*}
For $p>-C_0>-\infty$, with no loss of generality,
we can replace $p$ with $(p+C_0)$ throughout, since the Euler equation \eqref{e25}--\eqref{e26}
is invariant under the translation in $p$.
Thus, the analytic surface $(\mathcal{M}, \bg)$ corresponding to $(u,v,p)$ has positive curvature
for all $t \in [0,T]$.
On the other hand, in Lemma 4.2, we have another metric $\bg^\ast$ obtained from the shear flow.
Here $\bg^\ast$ is the induced metric of the following parameterised
map $\yy_{\bg^\ast}:\mathcal{M}' \subset \mathcal{M}\rightarrow \R^3$ near $\x^\star$:
\begin{equation*}
\yy_{\bg^\ast} = (Ax_2, Ax_1, f(x_2))^\top,
\end{equation*}
where $f$ is given implicitly by
$f'(x_2) = A \arctan \big(-A\int_0^{x_2}  u^2(s)\dd s\big)$ for $A > 1$.
We know that $\yy_\bg$ is strictly {\em short}:
	\begin{equation*}
	\yy_\bg^\#(\geucl) < (\yy_{\bg^\ast})^\#(\geucl)
	\end{equation*}
in the sense of quadratic forms, where $\#$ denotes the pullback operator.
Furthermore, as shown in Theorem 5.1, under condition \eqref{cauchy-k condition},
the geometric flow has an analytic solution, due to the Cauchy-Kowalewski theorem.
Thus, the previously constructed short map $\yy_\bg$ maps $[0,T]$
to $C^\infty_{\rm loc}(\mathcal{M}; \R^3)$.
	
\medskip
{\bf 2.} {\it Outline of the proof}.
Maps $(\Phi_1, \Phi_2, \Phi_3)$ will be explained in detail in the subsequent development;
they are given implicitly in \S 6.3 of \cite{CDS2012}, namely, in the proof of Corollary 1.2 therein.

\begin{lemma}[Corollary 1.2 in \cite{CDS2012}]
Let $n \in \mathbb{N}$ and let $\bg_0$ be a positive-definite $n \times n$ matrix.
There exists $r>0$ such that, for any smooth bounded open $\Omega \subset \R^n$
and any $\bg\in C^{0, \beta}(\overline{\Omega}; O^{+}(n))$ satisfying $\|\bg-\bg_0\|_{C^0} \leq r$,
the following holds{\rm :}

For any given $\yy_\bg \in C^1(\overline{\Omega}; \R^{n+1})$, $\e >0$,
and $\alpha \in (0, \min\{\frac{1}{1+2J_n}, \frac{\beta}{2}\})$,
there exists a map $\yy_{\rm w} \in C^{1,\alpha}(\overline{\Omega}; \R^{n+1})$ with
\begin{equation*}
\yy_{\rm w}^\# \geucl = \bg,\qquad
\|\yy_{\rm w} - \yy_\bg \|_{C^0} \leq \e.
\end{equation*}
\end{lemma}

Thus, in view of \eqref{arrows}, it is enough to prove that $(\Phi_1, \Phi_2, \Phi_3)$
are continuous in time, providing that all the function spaces therein
are endowed with the $C^1_{\rm loc}$ topology.
Here, as we begin with $\yy_\bg \in C^\infty$,
we are taking $n=2$, $\bg_0=\geucl$, and $\beta = 1$ in the lemma above.
Domain $\Omega$ is chosen to be a suitably small neighborhood
of $\x^\star$ in surface $\mathcal{M}$.

In the subsequent steps, we do not restrain ourselves to the case that $\mathcal{M}$
is a $2$-dimensional manifold immersed in $\R^3$.
Instead, the following arguments hold for any $n$-dimensional hypersurface $\mathcal{M}$
immersed into $(\R^{n+1}, \geucl)$, where $\geucl$ is the Euclidean metric.

\medskip
{\bf 3.} {\it Continuity of $\Phi_1$}.
In this step, we prove the continuity of $\Phi_1$, {\it i.e.}, the continuous dependence
of the Nash-Kuiper {\it wild} isometric immersions
with respect to the initial short immersion.
For simplicity of presentation, we only give the proof for Nash's construction in \cite{Nash1954},
which in fact requires at least two co-dimensions of the immersions.
Similar arguments work for Kuiper's construction in \cite{Kuiper} as well,
as long as we replace the ``Nash wrinkles'' (see \eqref{wrinkle} below)
by the ``corrugations'' in one co-dimension.
Our presentation of Nash's construction closely follows the exposition in \cite{exposition}.

Starting with the short map $\yy_\bg$, the $C^1$--isometric immersion is constructed
by adding ``Nash wrinkles'' to $\yy_\bg$ in countably many {\em stages},
and each stage involves finitely many {\em steps}:

\smallskip
\noindent
{\bf Nash's Steps.} To describe the steps, let us first recall the topological lemma concerning
the existence of a nice cover, which is proved by collecting the interiors of
the stars of (the barycentric subdivision of) a triangulation of $\mathcal{M}$:

\begin{lemma}[Lemma 2.2.1 in \cite{exposition}]\label{lemma on cover}
$\,$ Let $\mathcal{M}$ be an $n$-dimensional smooth manifold,
and $\{V_\lambda\}$ be an open cover.
Then there exists another cover $\{U_l\}$ such that
\begin{enumerate}
\item[\rm (i)]
Each $U_l$ lies in some $V_\lambda${\rm ;}
\item[\rm (ii)]
The closure of each $U_l$ is diffeomorphic to the $n$-dimensional closed ball in $\R^n${\rm ;}
\item[\rm (iii)]
Each $U_l$ intersects with at most finitely many other $U_{l'}$'s{\rm ;}
\item[\rm (iv)]
Each point $p\in\mathcal{M}$ has a neighbourhood contained in at most $(n+1)$ members of the cover{\rm ;}
\item[\rm (v)]
$\{U_l\}$ can be subdivided into $(n+1)$ classes, each consisting of pairwise disjoint $U_l$'s.
\end{enumerate}
\end{lemma}

Then, denoting by $I_l:=\{j: U_j \cap U_l \neq \emptyset\}$ in Lemma A(iii) with some $l$ fixed.
As $\yy_\bg$ is strictly short, for any $\delta>\|\bg-\yy_\bg^\# \geucl\|_{C^0}>0$,
we can choose $\delta_l>0$ such that
\begin{equation}\label{delta, step}
(1-\delta_l)\bg - \yy_\bg^\# \geucl \,\,\, \text{is positive-definite},\qquad\,\,
\|\delta_l \bg\|_{C^0(U_j)} \leq \frac{\delta}{2} \,\,\, \text{ for all } j \in I_l.
\end{equation}
Now, for some fixed $C^\infty$ partition of unity subordinate to $\{U_l\}$, we set
\begin{equation}\label{h-1}
\bh:=(1-\phi)\bg - \yy_\bg^\# \geucl,
\end{equation}
where $\phi := \sum_l \delta_l\phi_l$. Then $\bh$ is positive definite.
By Proposition 2.3.1 in \cite{exposition}, we can decompose $\bh$ into a locally finite
sum of {\em primitive metrics}:
\begin{equation}\label{decomp into primitive metrics}
\bh = \sum_j \bh_j
\end{equation}
such that each $\bh_j$ is supported in some $U_l$.
It is crucial to note that $\bh_j$ satisfies the following conditions:
For each $p\in\mathcal{M}$, there are at most $(n+1)J_n$ $\bh_j$'s
supported at $p$ and, for each $j$, ${\rm supp}(\bh_j)$ intersects
with finitely many other ${\rm supp}(\bh_{j'})$'s. Moreover,
\begin{equation}\label{primitive metric}
\bh_j = a_j^2 \dd\psi_j \otimes \dd\psi_j
\end{equation}
for some smooth functions $\phi_j$ depending only on $\mathcal{M}$
and the ``nice cover'' $\{U_l\}$ in Lemma \ref{lemma on cover}.

Now, choose two orthogonal vector fields $\boldsymbol{\nu}$
and $\boldsymbol{\xi}$ on $\overline{U_l}$,
which are of unit length and orthogonal to $TU_l$ throughout.
Then, in each {\em step} $j$ (in the sense of Nash), we add to $\yy_\bg$ a term $\ww^{\rm wrinkle}_{j}$,
which is a fast-oscillating plane wave of profile $a_j$,
frequency $\lambda\gg 1$ (to be determined), and directions $\boldsymbol{\xi}$
and $\boldsymbol{\nu}$.
More precisely, we consider
\begin{equation}\label{wrinkle}
\ww^{\rm wrinkle}_j(\x) = \frac{a_j(\x)}{\lambda}\cos(\lambda\psi_j(\x))\boldsymbol{\nu}(\x)
+ \frac{a_j(\x)}{\lambda} \sin (\lambda \psi_j (\x)) \boldsymbol{\xi}(\x).
\end{equation}
Such terms are known as ``Nash wrinkles'', or as ``spirals'' in Nash's original paper \cite{Nash1954}.

Finally, consider the map:
\begin{equation}\label{adding wrinkles}
\mathfrak{S} (\yy_\bg) := \yy_\bg + \ww^{\rm wrinkle}_1 + \ww^{\rm wrinkle}_2 + \ww^{\rm wrinkle}_3 + \cdots.
\end{equation}
To wit, for each point $\x\in \mathcal{M}$,
at most $(n+1)J_n$ Nash wrinkles are non-zero, so the sum in \eqref{adding wrinkles} is finite.
On the other hand, by choosing $\lambda$ sufficiently large,
we can require the Nash wrinkles to be very small in the $C^0$--norm.
Every such map $\mathfrak{S}$ is called a {\em Nash's stage}.

In summary, {\em $\mathfrak{S}$ maps from the space of $C^1_{\rm loc}$ strictly short
immersions $\mathcal{M} \hookrightarrow \R^3$ to itself}.

Now we are in the situation of showing that map $\mathfrak{S}$ is continuous in time.
In the sequel, $C_1, C_2, C_3,...$ denote universal constants depending only on the open cover
given by Lemma \ref{lemma on cover}.
Let us fix any $\eta > 0$ and assume that
\begin{equation}\label{x}
\|\yy_\bg(t)-\yy_\bg(s)\|_{C^1} \leq \eta \qquad \text{ for some } t,s \in [0,T].
\end{equation}
Then, in view of Eq. \eqref{h-1}, it follows that
\begin{equation*}
\|\bh_j(t)- \bh_j(s)\|_{C^1} \leq C_1 \eta^2 \qquad \text{ for all } j \in I_l.
\end{equation*}
Thanks to Eq. \eqref{primitive metric} (in which $\psi_j$ depends only on the cover), we have
\begin{equation*}
\|a_j(t)-a_j(s)\|_{C^1} \leq C_2 \eta \qquad \text{ for all } j \in I_l.
\end{equation*}
Therefore, expression \eqref{wrinkle} for the Nash wrinkles directly gives us
\begin{equation*}
\|\ww^{\rm wrinkle}_j (t) -\ww^{\rm wrinkle}_j (s)\|_{C^1} \leq \frac{2C_2}{\lambda} \eta \qquad \text{ for all } j \in I_l,
\end{equation*}	
which immediately implies that
\begin{equation}\label{xx}
\|\mathfrak{S}(\yy_{\bg}(t)) - \mathfrak{S}(\yy_{\bg}(s))\|_{C^1_{\rm loc}} \leq \frac{2(n+1)J_n C_2}{\lambda} \eta + \eta.
\end{equation}
From Eqs. \eqref{x}--\eqref{xx}, we conclude that map $\mathfrak{S}$ is continuous in time,
when its domain and range are equipped with the $C^1_{\rm loc}$ topology.
The above arguments hold for all $\lambda >0$; we are going to specify $\lambda$
in Eq. \eqref{choosing lambda} below, in order to ensure the convergence of the {\em stages}.

\smallskip
\noindent
{\bf Nash's Stages.}
The purpose of each stage $\mathfrak{S}$ is to correct the error,
$\|\mathfrak{S}(\yy_\bg)^\# \geucl - \bg\|_{C^0}$,
{\it i.e.} to lessen the deviation of the pulled back metrics from being isometric.
In view of Proposition 2.2.2 and the proof of Theorem 2.1.4 in \cite{exposition},
for any fixed $\e>0$,  we can obtain the following bound for $\mathfrak{S}$ at the $q$-th stage:
\begin{equation}\label{stage estimates}
\begin{cases}
\|\mathfrak{S}^q(\yy_\bg) - \mathfrak{S}^{q-1}(\yy_\bg) \|_{C^0(U_l)} < 2^{-q-1}\min\{\e, 2^{-l}\} \qquad \text{ for every } l, \\[1mm]
\|\bg-\mathfrak{S}^{q}(\yy_\bg)^\# \geucl\|_{C^0(\mathcal{M})} < \delta \equiv 4^{-q},\\[1mm]
\|D[\mathfrak{S}^q(\yy_\bg)] - D[\mathfrak{S}^{q-1}(\yy_\bg)]\|_{C^0(\mathcal{M})} < \sqrt{2}(n+1)J_n 2^{-q+1},
\end{cases}
\end{equation}
for each $q=1,2,3,\ldots$. Since we are proving everything locally, we assume without loss of generality
that $\mathcal{M}$ is compact.
By $\mathfrak{S}^q$ we mean the composition $\mathfrak{S}\circ\ldots\circ\mathfrak{S}$ for $q$ times.
As a remark, the above estimates involve the choice of $\lambda$ at the $q$-th stage
for each $q$ ({\it cf.} the proof of Eq. (2.16) in \cite{exposition}).

Hence, by the second inequality in Eq. \eqref{stage estimates}, we find that
\begin{equation*}
\Phi_1(\yy_\bg) := \lim_{q\rightarrow \infty} \mathfrak{S}^q(\yy_\bg) \in C^0(\mathcal{M}; \R^N)
\end{equation*}
is an isometric immersion.
Moreover, by the first and third equations, $\Phi_1(\yy_\bg)$ in fact lies in  $C^1(\mathcal{M}; \R^N)$.
All the above  constructions are {\it kinematic}, {\it i.e.}, they hold pointwise in $t \in [0,T]$.

Finally, in light of the proof for the estimates in \eqref{x}--\eqref{xx},
we have the following: If, at the $q$-th {\em stage},
we choose parameter $\lambda$ in the Nash wrinkles \eqref{wrinkle} to satisfy
\begin{equation}\label{choosing lambda}
\lambda = \lambda_q \geq 2^{q+1}(n+1)J_n C_2
\end{equation}
in addition to Eq. \eqref{stage estimates}, then we have
\begin{equation*}
\|\Phi_1(\yy_\bg)(t) - \Phi_1(\yy_\bg)(s)\|_{C^1_{\rm loc}} \leq 2\eta,
\end{equation*}
provided that $\|\yy_\bg(t)-\yy_\bg(s)\|_{C^1} \leq \eta$ for some  $t,s \in [0,T]$
as in Eq. \eqref{x}.

This completes the proof of the continuity of $\Phi_1$.

\smallskip
{\bf 4.} {\it Continuity of $\Phi_2$}.
This directly follows from the properties of mollification.
Indeed, let $0\leq J \in C^\infty(\R^n)$ be the standard mollifier in $\R^n$ such
that $\int_{\R^n}J(x)\,\dd x =1$ and ${\rm supp}(J) \in [-1,1]^n$. Then, for $\uu \in C^1({\Omega\subset \R^n; \R^N})$,
we define component-wise:
\begin{equation*}
\Phi_2(\uu)^i:= \uu^i \ast J_\e \qquad \text{in } \big\{\x\in\Omega: {\rm dist}(\x, \p\Omega) > \e\big\}\,\,
\text{ for } i\in\{1,2,\ldots,n\},
\end{equation*}
where $J_\e(\x):=\e^{-n}J(\frac{\x}{\e})$.
Then $\Phi_2(\uu)$ converges to $\uu$ in $C^1$ as $\e \rightarrow 0^{+}$ on any compact subset of $\Omega$.
In general, for $\uu \in C^1(\mathcal{M}; \R^N)$ where $\mathcal{M}$ is an $n$-dimensional manifold,
for any chart $\mathcal{M}' \subset \mathcal{M}$, we can find a $C^1$ diffeomorphism
$\ff: \mathcal{M}' \rightarrow \Omega \subset \R^n$.
Therefore, we have
\begin{equation*}
\|\Phi_2(\uu)(t) - \Phi_2(\uu)(s)\|_{C^1(\mathcal{M}')} \leq C\|\uu(t)-\uu(s)\|_{C^1(\mathcal{M})},
\end{equation*}
where $C$ depends on $\|\ff|_{\mathcal{M}'}\|_{C^1}$ and $\|(\ff|_{\mathcal{M}'})^{-1}\|_{C^1}$.
	
\smallskip
{\bf 5.} {\it Continuity of $\Phi_3$}. In this final step,
we prove that $\yy_{\rm w} = \Phi_3(\vv)$ is continuous.
This map is constructed by Conti-De Lellis-Sz\'{e}kelyhidi Jr. in \cite{{CDS2012}},
by a  {\em step/stage} construction  similar  to that of Nash's.
The difference is that, in every {\em step},
before adding the corrugations (introduced in Kuiper \cite{Kuiper} in co-dimension-1 case,
as the counterpart to the Nash wrinkles),
one first {\em mollifies}  immersion $\vv$ and metric $\bg$.
The estimates involved to control the mollification in each step is
motivated by Nash's argument for $C^\infty$ isometric embedding in \cite{nash2}.

\smallskip
\noindent
{\bf Steps in \cite{CDS2012}.}
Similar to Nash's construction, each {\em step} is achieved by adding the corrugations.
We closely follow \S 4 in \cite{CDS2012} for the presentation.
The basic building block is a corrugation function $\Gamma =\Gamma(z_1,z_2)\in C^\infty([0,\delta_\ast]\times\R;\R^2)$
which is $(2\pi)$--periodic in $z_2$ for some small $\delta_\ast >0$, and
\begin{equation}\label{estimate for Gamma}
\begin{cases}
|\p_{z_2}\Gamma(z_1,z_2) + (1,0)^\top|^2 = 1+z_1^2,\\[1mm]
|\p_{z_1}\p^k_{z_2}\Gamma_1(z_1, z_2)|+|\p_{z_2}^k\Gamma(z_1,z_2)| \leq C_k z_1 \qquad \text{ for } k \geq 0,
\end{cases}
\end{equation}
where $\Gamma=(\Gamma_1, \Gamma_2)^\top$.
As in Step 3 above, each {\em stage}, denoted by $\sss$ here,  consists of $J_n$ steps.
In each step, we add a {\em corrugation}:
\begin{equation*}
\sss (\vv):=\vv_0+\yy_1^{\rm corrugation} + \yy_2^{\rm corrugation} + \ldots + \yy_{J_n}^{\rm corrugation}.
\end{equation*}
Let us index the {\em steps} (in a fixed {\em stage}) by $j\in\{1,2,\ldots, J_n\}$, and abbreviate by
\begin{equation*}
\vv_j := \vv_0 + \yy_1^{\rm corrugation} + \ldots + \yy_{j}^{\rm corrugation},
\end{equation*}
where $\vv_0$ is given below.
Our goal is to describe $\vv_j$ or $\yy_j^{\rm corrugation}$ for each $j$, and investigate its dependence on time.

For this purpose, we first state the estimates achieved at the end of each {\em stage} in \cite{CDS2012}.
This will help us specifying the parameters ({\it i.e.}, $\lambda_j$, $l_j$, $a_j$, $\nu_j$, etc.) involved in each {\em step}:

\begin{lemma}[Proposition 5.1 in \cite{CDS2012}]\label{lemma: one stage, Phi3}
$\,$ For any $n\in\mathbb{N}$ and any positive-definite $n\times n$ matrix $\bg_0$,
there exists $r \in (0,1)$ such that, for any open bounded smooth $\Omega \subset \R^n$ and
any $\bg \in C^\beta(\overline{\Omega}; O^{+}(n))$ with $\|\bg-\bg_0\|_{C^0}\leq r$,
there is $\delta_0 > 0$ so that, for all $K \geq 1$, whenever
\begin{eqnarray*}
\|\vv^\# \geucl - \bg\|_{C^0} \leq \delta^2 \,\,\,\mbox{for some $\delta \leq \delta_0$},\qquad\,\,
\|\vv\|_{C^2} \leq \mu  \,\,\, \mbox{for some $\mu$},
\end{eqnarray*}
one can construct $\sss(\vv) \in C^2$ such that
\begin{eqnarray}\label{y1}
&&\|\sss(\vv)^\#\geucl - \bg\|_{C^0} \leq C_3 \delta^2\big(K^{-1} + \delta^{\beta-2}\mu^{-\beta}\big),\\\label{y2}
&&\|\sss(\vv)\|_{C^2} \leq C_3 \mu K^{J_n},\\ \label{y3}
&&\|\sss(\vv)-\vv\|_{C^1}\leq C_3 \delta,
\end{eqnarray}
where $C_3$ depends only on $n, \Omega, \bg_0$, and $\bg$.
\end{lemma}

It is crucial to remark that, in Lemma \ref{lemma: one stage, Phi3},
the construction of a certain {\em stage} begins with three given parameters:
$K, \delta$, and $\mu$.
In particular, they are independent of $j$, which indexes the {\em steps} within this {\em stage}.

From here, one first introduces the ``mollification parameter'' (Step 1 in \S 5.2, \cite{CDS2012}):
\begin{equation}\label{l, Phi3}
l:=\frac{\delta}{\mu},
\end{equation}
then, for the standard mollifying $0\leq J \in C^\infty_c(\R^n)$, we set
\begin{equation}\label{mollification, Phi3}
\vvv:=\vv\ast J_l, \qquad \gggg=\bg \ast J_l, \qquad J_l(\x):=l^{-n}J(\frac{\x}{l}).
\end{equation}
As $\vv$ is very close to being isometric (namely, $\|\vv^\#\geucl-\bg\|_{C^0} \leq \delta^2$),
for some large absolute constant $C_4$,  the following matrix $(1+\frac{C_4\delta^2}{r})\gggg - \vv^\#\geucl$
is positive-definite. Thus, it again can be decomposed into primitive metrics:
\begin{equation}\label{decomposition into primitive metrics, Phi3}
(1+\frac{C_4\delta^2}{r})\gggg - \vvv^\#\geucl = \sum_{i=1}^{J_n} \widetilde{a_i}^2\boldsymbol{\nu}_i \otimes \boldsymbol{\nu}_i.
\end{equation}
Then, we rescale
\begin{equation}\label{rescale, Phi3}
\vv_0 := \frac{1}{(1+C_4 r^{-1}\delta^2)^{1/2}}\vvv, \qquad a_i :=\frac{1}{(1+C_4 r^{-1}\delta^2)^{1/2}}\widetilde{a_i} \quad  \text{ for } i \in \{1,2,\ldots, J_n\}.
\end{equation}

Now we are ready for specifying each $\yy_j^{\rm corrugation}$ recursively.
The following is adapted from \S 4.2 of \cite{CDS2012}, by working a local orthonormal frame
$\{\mathbf{e}_1, \mathbf{e}_2, \ldots\}$ in $\R^n$.
First, define the vector fields:
\begin{equation}\label{two vector fields, Phi3}
\begin{cases}
\boldsymbol{\xi}_{j+1}:=\na \vv_{j} \cdot (\na^\top\vv_j\na\vv_j)^{-1}\cdot \boldsymbol{\nu}_j,\\
\boldsymbol{\zeta}_{j+1}:= \text{ the vector field dual to the $n$-form } \p_1 \vv_j \wedge \p_2\vv_j \wedge
\ldots \wedge \p_n  \vv_j.
\end{cases}
\end{equation}
Then the  ``amplitude'' is given by
\begin{equation}\label{amplitude, Phi3}
\Psi_j(x):=\frac{\boldsymbol{\xi}_j}{|\boldsymbol{\xi}_j|^2}(\x) \otimes \mathbf{e}_1 + \frac{\boldsymbol{\zeta}_j}{|\boldsymbol{\zeta}_j||\boldsymbol{\xi}_j|}(\x)\otimes \mathbf{e}_2.
\end{equation}
Finally, using the building block $\Gamma$, the {\em $j$-th corrugation} is defined by
\begin{equation}\label{corrugation, level j, Phi3}
\yy_j^{\rm corrugation} := \frac{1}{\lambda_j}\Psi_j(\x)\Gamma(|\boldsymbol{\xi}_j|a_j, \lambda_j \x\cdot \boldsymbol{\nu}_j),
\end{equation}
where, as in Step 3, \S 5.2 of \cite{CDS2012}, one chooses
\begin{equation}\label{lambda_j, Phi3}
\lambda_j := K^{j+1}l^{-1}.
\end{equation}

Let us now discuss the dependence of $\sss(\vv)$ on time.
For this purpose, fix $\eta >0$ and assume that
\begin{equation*}
\|\vv(t) - \vv(s)\|_{C^1} \leq \eta \qquad \text{ for some } t,s \in [0,T].
\end{equation*}
Then the mollification in equation \eqref{mollification, Phi3} gives us
\begin{align}\label{mollification estimate for v}
\|\vvv(t)-\vvv(s)\|_{C^2} &= \|\na J_l \ast \big(\vv(t)-\vv(s)\big)\|_{C^1} \nonumber\\
&= l^{-1} \Big\|\int_{\R^n}\na J (\z)\Big(\vv (t, \x-l\z) - \vv (s, \x-l\z)\Big)\,\dd \z \Big\|_{C^1}\nonumber\\
&\leq C_5 l^{-1} \|\vv(t)-\vv(s)\|_{C^1} \leq C_5 l^{-1}\eta,
\end{align}
where $C_5 \equiv \|J\|_{W^{1,1}(\R^n)}$.
From here,
the decomposition in Eq. \eqref{decomposition into primitive metrics, Phi3} yields
\begin{equation*}
\|\widetilde{a_i}(t) - \widetilde{a_i}(s)\|_{C^2} \leq C_6 \eta l^{-1},
\end{equation*}
where $C_6$ only depends on $\bg$, $\bg_0$, $n$, $\|J\|_{W^{1,1}(\R^n)}$, and the local geometry of $\Omega$.
Then \eqref{rescale, Phi3} shows that the rescaled quantities satisfy
\begin{equation}\label{v0, Phi3}
\|\vv_0(t)-\vv_0(s)\|_{C^2} + \|a_i(t)-a_i(s)\|_{C^2} \leq \frac{C_6 \eta}{l\sqrt{1+\frac{C_4\delta^2}{r}}}
\end{equation}
for each $i\in\{1,2,\ldots, J_n\}$. In addition, for the lower order derivatives, we have
\begin{equation}\label{lower order estimates for v, a. Phi3}
\|\vv_0(t)-\vv_0(s)\|_{C^1} + \|a_i(t)-a_i(s)\|_{C^1} \leq \frac{C_6 \eta}{\sqrt{1+\frac{C_4\delta^2}{r}}}.
\end{equation}

To proceed, notice that one can assume
\begin{equation}\label{assumption on nabla vj}
C_7^{-1}{\rm Id}\leq\na^\top\vv_j \na\vv_j\leq C_7{\rm Id}
\end{equation}
(see the beginning of \S 5.2 in \cite{CDS2012}); here $C_7$ may depend on $j$, but, as there are only finitely many $j$,
we can take $C_7$ to be absolute. Hence, equation \eqref{two vector fields, Phi3} implies that,
for each $j\in\{1,2,\ldots, J_n\}$,
\begin{equation}\label{estimate in xi and zeta}
\|\boldsymbol{\xi}_j(t)-\boldsymbol{\xi}_j(s)\|_{C^1} \leq  \frac{C_8\eta}{l\sqrt{1+\frac{C_4\delta^2}{r}}},
\quad \|\boldsymbol{\zeta}_j(t)-\boldsymbol{\zeta}_j(s)\|_{C^1}
\leq  \frac{C_8}{l}\bigg(\frac{\eta}{\sqrt{1+\frac{C_4\delta^2}{r}}}\bigg)^n,
\end{equation}
whereas the estimate in \eqref{lower order estimates for v, a. Phi3}
gives us
\begin{equation}\label{lower order estimate in xi and zeta}
\|\boldsymbol{\xi}_j(t)-\boldsymbol{\xi}_j(s)\|_{C^0} \leq  \frac{C_9\eta}{\sqrt{1+\frac{C_4\delta^2}{r}}},
\quad \|\boldsymbol{\zeta}_j(t)-\boldsymbol{\zeta}_j(s)\|_{C^0} \leq C_9 \bigg(\frac{\eta}{\sqrt{1+\frac{C_4\delta^2}{r}}}\bigg)^n.
\end{equation}

From here, we obtain the estimate for $\|\Psi_j(t)-\Psi_j(s)\|_{C^1}$: Since
\begin{equation*}
\na \frac{\boldsymbol{\xi}_j}{|\boldsymbol{\xi}_j|^2}
= \frac{\na \boldsymbol{\xi}_j - 2 \boldsymbol{\xi}_j \otimes \boldsymbol{\xi}_j}{|\boldsymbol{\xi}_j|^2},
\end{equation*}
and, by \eqref{assumption on nabla vj}, we have $|\boldsymbol{\xi}_j| \geq C_7^{-3/2}$ (and similarly for $\boldsymbol{\zeta}_j$)
so that
\begin{align*}
\|\Psi_j(t)-\Psi_j(s)\|_{C^1}
& \leq C_{10} \bigg\{\frac{\eta}{l\sqrt{1+\frac{C_4\delta^2}{r}}}
+\frac{1}{l}\bigg(\frac{\eta}{\sqrt{1+\frac{C_4\delta^2}{r}}}\bigg)^n\bigg\} \\
& \leq C_{10} \Big(\eta l^{-1} + \eta^n l^{-1}\Big).
\end{align*}

Next, let us bound $\|\yy_j^{\rm corrugation}(t)-\yy_j^{\rm corrugation}(s)\|_{C^1}$.
In view of expression \eqref{corrugation, level j, Phi3} of the corrugations,
a simple interpolation leads to
\begin{align*}
&\big\|\yy_j^{\rm corrugation}(t)-\yy_j^{\rm corrugation}(s)\big\|_{C^1} \nonumber\\
&\leq \frac{C_{11}}{\lambda_j} \bigg\{\big\|\Psi_j(t)-\Psi_j(s)\big\|_{C^1}\|(\widetilde\Gamma(t),\widetilde\Gamma(s))\|_{C^0}\\
&\qquad \quad\,\,\, + \big\| \widetilde\Gamma(t)-\widetilde\Gamma(s)\big\|_{C^1}\|(\Psi_j(t),\Psi_j(s))\|_{C^0} \bigg\},
\end{align*}
where the following shorthand is introduced:
$\widetilde\Gamma(s) \equiv \Gamma\big(|\boldsymbol{\xi}_j(s)|a_j(s), \lambda_j \x \cdot \boldsymbol{\nu}_j\big)$,
and similarly $\widetilde\Gamma(t)\equiv \Gamma\big(|\boldsymbol{\xi}_j(t)|a_j(t), \lambda_j \x \cdot \boldsymbol{\nu}_j\big)$.

To continue the estimate, we need the following uniform-in-time bounds.
First of all, thanks to \eqref{assumption on nabla vj}, we have
\begin{equation}\label{C0 estimate for xi, a, Psi}
\|\boldsymbol{\xi}_j\|_{C^0} + \|\Psi_j\|_{C^0} \leq C_{12},
\end{equation}
as well as
\begin{align*}
\|a_j\|^2_{C^1} &\leq \Big(1+C_4 r^{-1}\delta_2\Big)^{-1/2} \|\widetilde{a_j}\|_{C^1}
 \leq \|\widetilde{\bg} - \widetilde{\vv}^\sharp \bg_{\rm Eucl}\|_{C^1} + C_4 r^{-1} \delta^2 \|\widetilde{\bg}\|_{C^1}\nonumber\\
&\leq C_{5} l^{-1} \|\bg-\vv^\sharp \bg_{\rm Eucl}\|_{C^0} + C_4 C_5 r^{-1} \delta^2 \|\bg\|_{C^0},
\end{align*}
which can be proved in the similar manner to \eqref{mollification estimate for v}. It follows that
\begin{equation}\label{aj in C1 norm, Phi3}
\|a_j\|_{C^1} \leq C_{13} (1+l^{-1})\delta.
\end{equation}
Moreover, a simple computation gives
\begin{equation}\label{aj in C0 norm, Phi3}
\|a_j\|_{C^0} \leq C_{14} \delta.
\end{equation}
Finally, \eqref{assumption on nabla vj} yields that
\begin{equation*}
\|\boldsymbol{\xi}_j\|_{C^1} \leq C_{15} \Big(\|\vv_j\|_{C^0} + \|\vv_j\|_{C^1} \Big\|\na\Big[ \big(\na^\top \vv_j \na \vv_j\big)^{-1}\Big]\Big\|_{C^0}\bigg).
\end{equation*}
On the other hand, we have the identity:
\begin{align*}
&\na\Big[ \Big(\na^\top \vv_j \na \vv_j\Big)^{-1}\Big] \\
&= -\Big(\na^\top \vv_j \na \vv_j\Big)^{-1} \cdot \Big(\na\na^\top \vv_j \cdot \na \vv_j
+ \na^\top \vv_j \na^2\vv_j\Big) \cdot \Big(\na^\top \vv_j \na \vv_j\Big)^{-1},
\end{align*}
so that the following bound is verified:
\begin{equation}\label{xi C1 norm}
\|\boldsymbol{\xi}_j\|_{C^1} \leq C_{16} \|\vv_j\|_{C^2} \leq C_{17} \mu.
\end{equation}

Thus, as $\Gamma \in C^{\infty}([0,\delta_{\ast}]\times[0,2\pi])$ is periodic in the second arguments,
we have
\begin{align*}
&\big\|\yy_j^{\rm corrugation}(t)-\yy_j^{\rm corrugation}(s)\big\|_{C^1}\\
&\leq \frac{C_{18}}{\lambda_j} \Big\{ \big\|\Psi_j(t)-\Psi_j(s)\big\|_{C^1} + \big\|\widetilde\Gamma(t)-\widetilde\Gamma(s)\big\|_{C^1} \Big\} \nonumber\\
&\leq \frac{C_{19}}{\lambda_j} \Big\{\eta l^{-1} + \eta^nl^{-1} + \big\|\widetilde\Gamma(t)-\widetilde\Gamma(s)\big\|_{C^1} \Big\}.
\end{align*}
To continue, let us estimate by the Taylor expansion
\begin{align*}
\big\|\widetilde\Gamma(t)-\widetilde\Gamma(s)\big\|_{C^1}
\leq \Big\| \p_{z_1}\Gamma(\Theta, \lambda_j \x \cdot \boldsymbol{\nu}_j) \Big[|\boldsymbol{\xi}_j(t)|a_j(t) - |\boldsymbol{\xi}_j(s)|a_j(s)\Big] \Big\|_{C^1},
\end{align*}
where $|\Theta|$ lies between $|\boldsymbol{\xi}_j(s)||a_j(s)|$ and $|\boldsymbol{\xi}_j(t)||a_j(t)|$. Then
\begin{align*}
&\big\|\widetilde\Gamma(t)-\widetilde\Gamma(s)\big\|_{C^1}\\
&\leq C_{20}\Big\||\boldsymbol{\xi}_j(t)|a_j(t) - |\boldsymbol{\xi}_j(s)|a_j(s) \Big\|_{C^1} \nonumber\\
&\leq C_{21} \bigg\{\|\boldsymbol{\xi}_j(t)\|_{C^1} \|a_j(t)-a_j(s)\|_{C^0} + \|\boldsymbol{\xi}_j(t)\|_{C^0} \|a_j(t)-a_j(s)\|_{C^1}\nonumber\\
&\qquad\quad\,\,  + \|a_j(s)\|_{C^1} \Big\| |\boldsymbol{\xi}_j(t)|-|\boldsymbol{\xi}_j(s)|\Big\|_{C^0}
+ \|a_j(s)\|_{C^0}\Big\| |\boldsymbol{\xi}_j(t)|-|\boldsymbol{\xi}_j(s)|\Big\|_{C^1} \bigg\}\nonumber\\
&\leq C_{22} \Big\{\mu \delta + \eta + (1+l^{-1})\delta\eta + \delta\eta l^{-1} \Big\}.
\end{align*}
Here, the first line follows from \eqref{estimate for Gamma}, \eqref{C0 estimate for xi, a, Psi},
and \eqref{aj in C0 norm, Phi3},
the third inequality follows from interpolation, and the final one from
estimates \eqref{estimate in xi and zeta}--\eqref{xi C1 norm}.
The constants $C_{20}, C_{21}$, and $C_{22}$ may further depend on $\|\p_{z_1}\Gamma\|_{C^0}$.
We are now ready to conclude
\begin{equation}\label{corrugation difference, Phi3}
\big\|\yy_j^{\rm corrugation}(t)-\yy_j^{\rm corrugation}(s)\big\|_{C^1}
\leq \frac{C_{23}}{\lambda_j}\big\{ \eta l^{-1} + \mu \delta +  \eta + l^{-1}\delta \eta\big\},
\end{equation}
where, without loss of generalities, we have assumed that $\eta \leq 1$.

Therefore, summing over the geometric series in view of equation \eqref{lambda_j, Phi3}, we have
\begin{align}
\|\sss(\vv)(t)-\sss(\vv)(s)\|_{C^1}
&\leq C_{23} \frac{1-K^{-J_n}}{K^2-K} (\eta + \mu l \delta + \eta l + \delta\eta)\nonumber\\
&\leq C_{23} (\eta + \mu l \delta + \eta l + \delta\eta),\label{one stage, Phi3}
\end{align}
since $K>1$ in the assumption of Lemma \ref{lemma: one stage, Phi3}.
Here, $C_{23}$ depends on $\bg_0$, $\bg$, $r$, $n$, $\Omega$, $\beta$,  $\|J\|_{W^{1,1}}$, and $\Gamma$,
but not on $\eta$, $l$, $\mu$, and $\delta$.
The last three of these four parameters are chosen differently for distinctive {\em stages} below.

\medskip
\noindent
{\bf Stages in \cite{CDS2012}.}	
Now we iterate for countably many times the construction above for one {\em stage}.
As in Step 3, we index the stages by $q=1,2,3,\ldots$.
Recall from \S 6 in \cite{CDS2012} that $\Phi_3$ is given by
\begin{equation*}
\Phi_3(\vv):= \lim_{q\rightarrow \infty} \sss^q (\vv) \in C^1_{\rm loc}(\mathcal{M}; \R^{n+1}),
\end{equation*}
where one needs to suitably choose $(\mu_q, \delta_q)$ in place of $(\mu, \delta)$
in equations \eqref{y1}--\eqref{y3} in Lemma \ref{lemma: one stage, Phi3} above,
with $\sss$ is replaced by $\sss^q=\underbrace{\sss\circ\ldots\circ\sss}_{q \text{ times}}$ therein.

Following the delicate arguments therein, by choosing
\begin{equation*}
a < \min \big\{\frac{1}{2}, \frac{\beta J_n}{2-\beta} \big\}, \qquad \alpha < \min \big\{\frac{\beta}{2}, \frac{1}{1+2J_n}\big\},
\end{equation*}
one can bound via the interpolation of the $C^1$ and $C^2$ estimates as follows: For each $q$,
\begin{align}
& \|\sss^{q+1}(\vv) - \sss^q(\vv)\|_{C^{1,\alpha}} \leq C_{24} \mu_0 K^{-\big((1-\alpha)a-\alpha J_n\big)q} \quad \text{ on } [0,T],
         \label{holder}\\
&\|\sss^q(\vv)^\# \geucl - \bg\|_{C^0} \leq \delta_q^2.\label{isometric for holder case}
\end{align}
Therefore, in view of equation \eqref{holder},
$\Phi_3(\vv)$ in fact lies in $C^{1,\alpha}_{\rm loc}(\mathcal{M}; \R^{n+1})$
and, by Eq. \eqref{isometric for holder case}, $\Phi_3(\vv)$ is indeed an isometric immersion.
		
It remains to discuss the dependence of $\Phi_3(\vv)$ on time.
From the definition of $\Phi_3$ and estimate \eqref{one stage, Phi3} for one stage,
we observe
\begin{equation}\label{Phi 3 bound in t}
\|\Phi_3(\vv)(t)-\Phi_3(\vv)(s)\|_{C^1} \leq \lim_{q\rightarrow \infty} C_{23} (\eta + \mu_q l_q \delta_q + \eta l_q + \eta \delta_q).
\end{equation}
In \S 6.1 of \cite{CDS2012}, parameters $(\delta_q, \mu_q)$ are chosen to satisfy
\begin{equation}\label{relation for delta, mu}
\delta_{q} \leq \delta_0 K^{-aq}, \qquad \mu_q = \mu_0 K^{q J_n},
\end{equation}
where $\mu_0 > 0$ is a fixed constant, and $K \geq 2^{1/a}$ for $a>0$ specified in the preceding.
Since $l_q:=\delta_q \mu^{-1}_q$ as defined in equation \eqref{l, Phi3}, we have
\begin{equation*}
l_q \leq \frac{\delta_0}{\mu_0}K^{-(a+J_n)q} \qquad \text{ for each } q\in\{1,2,3\ldots\}.
\end{equation*}
In particular, $l_q \rightarrow 0$ because $K>1$. Therefore, \eqref{Phi 3 bound in t} becomes
\begin{equation*}
\|\Phi_3(\vv)(t)-\Phi_3(\vv)(s)\|_{C^1} \leq C_{23} \eta,
\end{equation*}
where $C_{23}$ is a universal constant depending on $\bg_0$, $\bg$, $n$, $r$, $\Omega$, $\delta_0$, $\mu_0$, $K$, $\alpha$, $\beta$,
$\Gamma$, $\|J\|_{W^{1,1}}$, and $a$,
provided that $\|\vv(t) - \vv(s)\|_{C^1} \leq \eta$ for some $t,s \in [0,T]$.
That is, $\Phi_3: C^2_{\rm loc}(\mathcal{M}; \R^{n+1}) \rightarrow C^{1,\alpha}_{\rm loc}(\mathcal{M}; \R^{n+1})$
is continuous in time, when the domain and range are both equipped with the $C^1_{\rm loc}$ topology.
The proof is now complete.	
\endproof
\end{appendices}

\bigskip

\section*{Acknowledgments}
A. Acharya acknowledges the support of the Rosi and Max Varon Visiting Professorship at
the Weizmann Institute of Science, Rehovot, Israel,
and was also supported in part by grants NSF-CMMI-1435624, NSF-DMS-1434734,
and ARO W911NF-15-1-0239.
G.-Q. Chen's research was supported in part by
the UK
Engineering and Physical Sciences Research Council Award
EP/E035027/1 and
EP/L015811/1, and the Royal Society--Wolfson Research Merit Award (UK).
S. Li's research was supported in part by the UK EPSRC Science and
Innovation award to the Oxford Centre for Nonlinear PDE (EP/E035027/1).
M. Slemrod was supported in part by Simons Collaborative Research Grant 232531.
M. Slemrod also thanks the Oxford Center for Nonlinear PDE and the Max Planck Institute for Mathematics  in the Sciences (Leipzig) for their kind hospitality.
D. Wang was supported in part by NSF grants DMS-1312800 and DMS-1613213.
We also thank L. Sz\'{e}kelyhidi Jr.  for his valuable remarks and suggestions.
Finally we thank the anonymous referee for his/her valuable comments and suggestions.

\end{document}